\documentclass[10pt]{amsart}
\usepackage{amssymb, amsthm, amsmath}
\usepackage{latexsym}

\theoremstyle{plain}
\newtheorem{prop}{Proposition}
\newtheorem{thm}{Theorem}
\newtheorem{cor}{Corollary}
\newtheorem{lemma}{Lemma}
\newtheorem*{movlemma}{Moving Lemma}

\theoremstyle{remark}

\theoremstyle{definition}
\newtheorem{defn}{Definition}

\def\bP{{\mathbb P}}

\def\Z{{\mathbb Z}}

\def\R{{\mathbb R}}
\def\C{{\mathbb C}}
\def\cO{{\mathcal O}}
\def\F{{\mathcal F}}

\def\E{{\mathcal E}}

\def\D{{\mathcal D}}

\def\CS{{\mathfrak C}}

\def\QQ{{\mathcal Q}}
\def\ev{{\mathrm{ev}}}
\def\WW{{\widehat W}}
\def\l{{\lambda}}
\def\m{{\mu}}
\def\n{{\nu}}
\def\L{{\Lambda}}

\def\X{{\mathfrak X}}

\def\a{{\alpha}}

\def\s{{\sigma}}
\def\t{{\tau}}

\DeclareMathOperator{\Spec}{Spec}
\DeclareMathOperator{\rk}{rk}

\DeclareMathOperator{\Span}{Span}

\newcommand{\dis}{\displaystyle}

\newcommand{\ssm}{\smallsetminus}
\newcommand{\gequ}{\geqslant}
\newcommand{\lequ}{\leqslant}
\newcommand{\lra}{\longrightarrow}

\newcommand{\ra}{\rightarrow}

\newcommand{\ov}{\overline}
\newcommand{\noin}{\noindent}
\newcommand{\wt}{\widetilde}
\newcommand{\Pf}{\mbox{Pfaffian}}

\begin{document}

\title[Quantum Cohomology of the Lagrangian Grassmannian]
{Quantum Cohomology of the Lagrangian Grassmannian}
\author{Andrew Kresch and Harry Tamvakis}
\date{May 3, 2003}
\address{Department of Mathematics, University of Pennsylvania,
209 South 33rd Street,
Philadelphia, PA 19104-6395}
\email{kresch@math.upenn.edu, harryt@math.upenn.edu}

\begin{abstract}
Let $V$ be a symplectic vector space and $LG$ be the Lagrangian
Grassmannian which  parametrizes maximal isotropic
subspaces in $V$. We give a presentation for the (small) quantum
cohomology ring $QH^*(LG)$ and show that
its multiplicative structure is determined by the ring of 
$\wt{Q}$-polynomials. We formulate a 
`quantum Schubert calculus' which includes quantum Pieri
and Giambelli formulas, as well as algorithms for computing the
structure constants appearing in the quantum product of 
Schubert classes.
\end{abstract}

\maketitle

\section{Introduction}
\label{intro}

\noindent
The multiplicative structure
of the quantum cohomology ring $QH^*(\X)$ of a projective complex 
manifold $\X$ encodes the enumerative
geometry of rational curves in
$\X$, in the form of {\em Gromov--Witten invariants}. The ring $QH^*(\X)$
is a deformation of the cohomology ring $H^*(\X)$ which first
appeared in the work of string theorists (see e.g. \cite{V}, \cite{W}).
The exposition in \cite{FPa},
following the work of Kontsevich and Manin \cite{KM}, gives an
algebro-geometric approach to 
this theory when $\X$ is a homogeneous space $G/P$, where $G$ is a
complex Lie group and $P$ a parabolic subgroup, which is
our main interest here (see also \cite{LT}).
We will work throughout with the small
quantum cohomology ring (terminology from \cite[\S 10]{FPa}). 

The cohomology rings of homogeneous spaces $\X=G/P$ have been studied
extensively; see \cite{Bo} \cite{BGG} \cite{D1} \cite{D2}. We would
like to have an
analogous description of the multiplicative structure of $QH^*(\X)$,
which we refer to as {\em quantum Schubert calculus}. As in the 
classical case, there are three main ingredients necessary for
the latter theory: (i) a {\em presentation} of 
$QH^*(\X)$ in terms of generators and relations, (ii) 
a {\em quantum Giambelli
formula} which identifies the polynomials which represent the Schubert 
classes in this presentation, and (iii) algorithms for computing the
structure constants in the multiplication table of $QH^*(\X)$ (the 
latter include {\em quantum Pieri rules}). 

Currently we have a fairly complete
understanding these questions when $\X=SL_n/P$ is a partial flag
variety of $SL_n(\C)$. Since the theory is not functorial, a separate
analysis must be done for the various parabolic subgroups $P$.
For work (in the $SL_n$ case) on the quantum cohomology of (i)
Grassmannians, see \cite{W}, \cite{ST}, \cite{Ber}, \cite{BCF},
(ii) complete flag varieties, see \cite{GK}, \cite{CF1},
\cite{FGP} and (iii) partial flag varieties, see \cite{AS}, \cite{K1},
\cite{K2}, \cite{CF2}, \cite{C}. 

In contrast to the $SL_n$ situation, much less is understood
for the other families 
of Lie groups. For an arbitrary complex semisimple Lie group $G$,
with Borel subgroup $B$,
Kim \cite{K3} found the
quantum $\mathcal D$-module structure for the flag variety
$\X=G/B$, and thus determined a presentation of $QH^*(\X)$.
To the authors' knowledge one still lacks a presentation of the quantum
ring for the other parabolic subgroups and a quantum Giambelli
formula to compute the Gromov--Witten invariants, if $G$ is
not of type $A$. 
Our aim in this paper is to answer all three
questions when the Lie group $G=Sp_{2n}$ is the symplectic group
and $P=P_n$ is the maximal 
parabolic subgroup associated with a `right end root' in the Dynkin
diagram for the root system of type $C_n$. 
In this case $\X=G/P$ is the Lagrangian 
Grassmannian described below. In a companion paper to this one 
\cite{KTorth}, we describe
the corresponding story for the orthogonal groups.

All cohomology classes in this paper occur in even degrees.
To avoid unnecessary factors of two, we adopt the following convention:
a class $\a$ in the cohomology of a complex variety $\X$
has degree $k$ when $\a$ lies in
$H^{2k}(\X,\Z)$.

Let $V$ be a complex vector space of dimension $2n$, equipped with a 
symplectic form. The variety of Lagrangian (i.e., maximal isotropic)
subspaces of $V$ is 
the {\em Lagrangian Grassmannian}
$LG=LG(n,2n)=Sp_{2n}/P_n$. 
The integral cohomology ring $H^*(LG,\Z)$ 
has a $\Z$-basis of Schubert classes $\s_{\l}$, 
one for each strict partition $\l=(\l_1>\l_2>\cdots>\l_r)$ with
$\l_1\lequ n$. Each $\s_{\l}$ 
is the class of a Schubert variety $\X_{\l}$ 
of codimension $|\l|=\sum\l_i$ in $LG$. 
Of particular importance are the classes
$\s_{i,j}$, 
determined by two Schubert conditions, and the
special Schubert classes $\s_i=\s_{i,0}$,
defined by a single
Schubert condition. This terminology is reviewed in \S
\ref{schuvar}.

To describe the multiplicative structure of the classical and quantum 
cohomology ring of $LG$ we will use the 
{\em $\wt{Q}$-polynomials} of Pragacz and Ratajski
\cite{PR}. These symmetric polynomials are modeled on Schur's
$Q$-functions \cite{S}, and are the geometric analogues
of Schur's $S$-functions in type $C$.
Let $X_k=(x_1,\ldots,x_k)$ be a $k$-tuple of variables, let
$X=X_n$, and
for $i\gequ 0$ let $\wt{Q}_i(X)$ be the 
$i$-th elementary symmetric polynomial $e_i(X)$.
For any nonnegative integers  $i,j$  with $i\gequ j$, let
\begin{equation}
\label{clgiam1}
\wt{Q}_{i,j}(X)=
\wt{Q}_i(X)\wt{Q}_j(X)+
2\sum_{k=1}^j(-1)^k\wt{Q}_{i+k}(X)\wt{Q}_{j-k}(X),
\end{equation}
and for any partition
$\l=(\l_1\gequ\l_2\gequ\cdots\gequ\l_r\gequ 0)$, define
\begin{equation}
\label{clgiam2}
\wt{Q}_{\l}(X)=\Pf[\wt{Q}_{\l_i,\l_j}(X)]_{1\lequ i<j\lequ r},
\end{equation}
if $r$ is even; we may always assume this is the case by setting
$\l_r=0$ if necessary. For each positive
integer $n$ we denote by $\D_n$ the set of strict
partitions $\l$ with $\l_1\lequ n$.

Let $\L_n$ denote the ring $\Z[X]^{S_n}$
of symmetric polynomials in $X$; as an abelian group, $\L_n$ 
is spanned by the polynomials $\wt{Q}_{\l}(X)$ for all 
partitions $\l$. From the analysis
in \cite[Sect.\ 6]{P} and \cite{PR} we obtain that the map
which sends $\wt{Q}_{\l}(X)$
to $\s_{\l}$ for all $\l\in\D_n$ extends to a 
surjective {\em ring} homomorphism $\phi:\L_n\ra H^*(LG,\Z)$ 
with kernel generated by the relations $\wt{Q}_{i,i}(X)=0$ 
for $1\lequ i\lequ n$; the morphism $\phi$ is evaluation on 
the Chern roots of the tautological rank $n$ quotient 
bundle over $LG$. From this
statement one gets a presentation for 
the cohomology ring of $LG$; moreover,  
equations (\ref{clgiam1}) and (\ref{clgiam2})
become Giambelli-type formulas, expressing the Schubert classes in terms
of the special ones.

The quantum cohomology $QH^*(LG)$ of
the Lagrangian Grassmannian $LG(n,2n)$
is an algebra over $\Z[q]$, where $q$ is a formal
variable of degree $n+1$; one recovers 
the classical cohomology ring by setting $q=0$. 
For the $SL_n$-Grassmannian, Bertram \cite{Ber} proved the
remarkable fact that the classical and quantum Giambelli
are identical; note that this is in contrast to the  other 
$SL_n$-flag varieties (\cite{FGP}, \cite{CF2}, \cite{C}).
Strictly speaking, 
the quantum Giambelli formulas for $LG$ do not coincide with the 
classical ones. Let $X^+:=X_{n+1}$, and define $\wt{\L}_{n+1}$
to be the subring of $\L_{n+1}$ generated by the polynomials
$\wt{Q}_i(X^+)$ for $i\lequ n$ together with the polynomial 
$2\,\wt{Q}_{n+1}(X^+)$; the latter will play the role of $q$
in the quantum cohomology ring.
Equations (\ref{clgiam1}) and (\ref{clgiam2}) imply that
$\wt{Q}_{\l}(X^+)\in\wt{\L}_{n+1}$ for all partitions $\l$ 
with $\l_1\lequ n$.

\begin{thm}
\label{lgthm} 
The map sending
$\wt{Q}_{\l}(X^+)$ to $\s_{\l}$ for all $\l\in\D_n$
and $2\,\wt{Q}_{n+1}(X^+)$ to $q$ extends to 
a surjective {\em ring} homomorphism $\wt{\L}_{n+1}\ra QH^*(LG)$
with kernel 
generated by the relations $\wt{Q}_{i,i}(X^+)=0$ for $1\lequ i\lequ n$.
The ring $QH^*(LG)$ is presented as a quotient of the polynomial
ring $\Z[\s_1,\ldots,\s_n,q]$ by the relations 
\begin{equation}
\label{lgqrel}
\s_i^2+2\sum_{k=1}^{n-i}(-1)^k\s_{i+k}\s_{i-k}=(-1)^{n-i}\s_{2i-n-1}\,q,
\end{equation}
where it is understood that $\s_j=0$ for $j<0$.
The Schubert class $\s_{\l}$
in this presentation is given by the Giambelli formulas
\begin{equation}
\label{lgqgiam1}
\s_{i,j}=
\s_i\s_j+2\sum_{k=1}^{n-i}(-1)^k\s_{i+k}\s_{j-k}+(-1)^{n+1-i}\s_{i+j-n-1}\,q
\end{equation}
for $i>j>0$, and 
\begin{equation}
\label{lgqgiam2}
\s_{\l}=\text{\em Pfaffian}[\s_{\l_i,\l_j}]_{1\lequ i<j\lequ r},
\end{equation}
where quantum multiplication is employed
throughout. In other words, quantum Giambelli for $LG(n,2n)$ 
coincides with classical Giambelli for $LG(n+1,2n+2)$, when the
class $2\,\s_{n+1}$ is identified with $q$.
\end{thm}

Our proof of Theorem \ref{lgthm}
is similar to \cite{Ber} in that we use the
Lagrangian Quot scheme $LQ_d$ to compactify 
the moduli space of degree $d$ maps $\bP^1\ra LG$. 
However, unlike the Quot scheme in type
$A$, the scheme $LQ_d$ is singular in general. More
significantly, in
\cite{KT} we showed that there is no direct analogue of the 
Kempf--Laksov formula \cite{KL} for {\em isotropic morphisms}
in type $C$. Thus the key intersection-theoretic ingredient
used by Bertram \cite{Ber} to prove quantum Giambelli in 
type $A$ is no longer available. Instead, we first show how to
evaluate any three term Gromov--Witten invariant in degree $1$
(i.e., one that counts lines on $LG$ satisfying three incidence
conditions); this
provides the quantum Giambelli formula (\ref{lgqgiam1}) for
the Schubert classes $\s_{i,j}$, as well as the
ring presentation for $QH^*(LG)$ (\S \ref{countinglines}).
The Pfaffian formula (\ref{lgqgiam2})
then follows from a relation in the Chow group of Lagrangian 
Quot schemes, which is proved using a careful analysis of the
boundary of $LQ_d$ and some remarkable Pfaffian identities 
for certain {\em symplectic Schubert polynomials} (\S
\ref{pfids}). The latter objects were introduced in \cite{PR}
and represent the Schubert classes in the 
complete symplectic flag variety.

In $QH^*(LG)$ there are formulas
\[
\s_{\l}\cdot \s_{\m} = 
\sum e_{\l\m}^{\n}(n)\,\s_{\n}\, q^d,
\]
the sum over $d\gequ 0$ and strict partitions $\n$ with 
$|\n|=|\l|+|\m|-d(n+1)$. The {\em quantum structure constant}
$e_{\l\m}^{\n}(n)$ is equal to a Gromov--Witten invariant
$\langle \s_{\l}, \s_{\m}, \s_{\n'} \rangle_d$, which is a
non-negative integer. The integer
$\langle \s_{\l}, \s_{\m}, \s_{\n'} \rangle_d$ counts the number
of rational curves of degree $d$ that meet three Schubert varieties
$\X_{\l}$, $\X_{\m}$ and $\X_{\n'}$ in general position. 
Here $\n'$ is the dual partition
of $\n$, defined so that the classes $\s_{\n}$ and $\s_{\n'}$
are Poincar\'e dual to each other in $H^*(LG,\Z)$. 

Using Theorem \ref{lgthm}, we give
formulas and algorithms to compute the quantum numbers
$e_{\l\m}^{\n}(n)$; our approach was
inspired by the second author's study \cite{T2} of Arakelov theory
on Lagrangian Grassmannians. The resulting `quantum Schubert calculus'
is a natural extension of the classical formulas of Hiller and Boe
\cite{HB}, Stembridge \cite{St} and Pragacz \cite{P}. 
In particular we obtain a 
`quantum Pieri rule' (Proposition \ref{qupieri}) for
multiplying by a special Schubert class in $QH^*(LG)$, and some
cases of a general `quantum Littlewood-Richardson rule'
in type $C$, which we obtain by relating some of the higher degree
quantum structure constants with classical ones (Theorem 
\ref{symmetrythm} and Corollary \ref{QLRcor}).

This paper is organized as follows. In Section \ref{Qtilde} we study the 
$\wt{Q}$-polynomials
 and symplectic Schubert polynomials, and prove the Pfaffian
identities we need for the latter family. The Lagrangian and isotropic
Grassmannians are introduced in Section \ref{isograss}, which includes our
proof of the easier part (\ref{lgqgiam1}) of the quantum Giambelli
formula for $LG$. We complete the proof of Theorem \ref{lgthm} in
Sections \ref{lagrquot} and  \ref{itlqd}, by studying intersections
on the Lagrangian Quot scheme $LQ_d$. Finally, in Section \ref{qsc} we
formulate a `quantum Schubert calculus' for $LG$ which extends the
classical one.

The authors thank Anders Buch and Bill Fulton for helpful correspondence.
Both authors were supported in part by National Science
Foundation post-doctoral research fellowships.
We also wish to
thank the Institut des Hautes \'Etudes Scientifiques
for its hospitality during the final stages of this work. 

\section{$\wt{Q}$-polynomials and symplectic Schubert polynomials}
\label{Qtilde}

\subsection{Basic definitions and properties}
\label{bdefs}
We begin by recalling some basic facts about partitions and their
Young diagrams; the main reference is \cite{M}. 
A {\em partition} is a sequence $\l=(\l_1,\l_2,\ldots,\l_r)$
of nonnegative integers in decreasing order; the {\em length} 
$\ell(\l)$ is the number of (nonzero) parts $\l_i$,
and the {\em weight} $|\l|=\sum\l_i$. Set $\l_k=0$ for any $k>\ell(\l)$.
We identify a partition with its associated {\em Young diagram} of boxes;
the relation $\l\supset\mu$ is defined by the containment of
diagrams. When $\l\supset\mu$
the set-theoretic difference is the {\em skew
diagram} $\l/\m$. 
A skew diagram $\alpha$ is a {\em horizontal strip} if it has at most
one box in each column. Two boxes in $\alpha$ are {\em connected}
if they share a vertex or an edge; 
this defines the connected components of $\alpha$.

Let $\E_n$ denote the set of all partitions $\l$ with
$\l_1\lequ n$.
A partition is {\em strict} if all its parts are different.
Let $\rho_n=(n,n-1,\ldots,1)$ and recall that $\D_n$ denotes
the set of strict partitions $\l$ with $\l\subset\rho_n$.
We use the notation $\l\ssm\m$ to denote the partition with parts given
by the parts of $\l$ which are not parts of $\m$.


Our references for the polynomials in this section are
\cite{PR} and \cite{LP1}; we will follow the notational 
conventions of \cite{KT} throughout most of the paper. 
Recall from the introduction the definitions of the polynomials
$\wt{Q}_{\l}(X)$ and the ring
$\L_n$ of symmetric polynomials in $X=X_n$. These polynomials
enjoy the following properties \cite[Sect.\ 4]{PR}:

\medskip
(a) If $\l_1>n$, then $\wt{Q}_{\l}(X)=0$.

\medskip
(b) The set $\{\wt{Q}_{\l}(X)\ |\ 
\l\in\E_n\}$ is a $\Z$-basis for $\L_n$.

\medskip
(c) $\wt{Q}_{i,i}(X)=e_i(X_1^2,\ldots,X_n^2)$ for all $i$.

\medskip
(d) $\wt{Q}_n(X)\wt{Q}_{\l}(X)=\wt{Q}_{(n,\l)}(X)$ for all $\l\in\E_n$.

\medskip
(e) If $\l=(\l_1,\ldots,\l_r)$ and 
$\l^+=\l\cup(i,i)=(\l_1,\ldots,i,i,\ldots,\l_r)$ then 
\[
\wt{Q}_{\l^+}(X)=\wt{Q}_{i,i}(X)\wt{Q}_{\l}(X).
\]

Define a {\em composition} $\nu$ to be 
a sequence of nonnegative integers with finitely many nonzero parts;
we set $|\nu|=\sum \nu_i$. The $\wt{Q}$-polynomials
make sense when indexed by a composition $\nu$
(in fact, the index can be an arbitrary
finite-length sequence of integers, following
\cite[Rmk.\ 5.3]{LP1} --- we will not need this here).
For any composition $\nu$
the polynomial $\wt{Q}_{\nu}(X)$ is equal to $\pm\wt{Q}_{\l}(X)$, where 
$\l$ is the partition obtained from $\nu$ by reordering the parts of
$\nu$. The sign is determined by the following rule:
\[
\wt{Q}_{(\ldots,i,j,\ldots)}(X)=-\wt{Q}_{(\ldots,j,i,\ldots)}(X)
\]
for any two {\em adjacent} indices $i$, $j$ with $i\neq j$. In other
words, for any composition $\nu$ we have $\wt{Q}_{\nu}(X)=
\mathrm{sgn}(\s)\wt{Q}_{\l}(X)$, where $\s\in S_{\infty}$ is the
permutation of minimal length such that $\nu_{\s(k)}=\l_{k}$ for
all $k$. The polynomials defined in this way 
satisfy Pfaffian relations such as 
\begin{equation}
\label{pfaff}
\wt{Q}_{\nu}(X)=
\sum_{j=1}^{r-1}(-1)^{j-1}\wt{Q}_{\nu_j,\nu_r}(X)
\cdot \wt{Q}_{\nu\ssm\{\nu_j,\nu_r\}}(X),
\end{equation}
where $r$ is an even number such that $\nu_i=0$ for $i>r$;
this follows e.g.\ from
the Laplace-type expansion for Pfaffians displayed in \cite[(D.1)]{FP}.
Using this convention, we have the following extension of 
\cite[Prop.\ 4.1]{PR}:

\begin{prop}
\label{ex4.1}
For any partition $\l$ (not necessarily strict) we have
\begin{equation}
\label{extn}
\wt{Q}_{\l}(X)=
\sum_{k=0}^{\ell(\l)}x_1^k\sum_{\mu}\wt{Q}_{\mu}(X'),
\end{equation}
where the inner sum is over all compositions $\mu$ such that 
$|\l|-|\mu|=k$ and $\l_i-\mu_i\in\{0,1\}$ for each $i$,
while $X'=(x_2,\ldots,x_n)$.
\end{prop}
\begin{proof}
One notices the proof of \cite[Prop.\ 4.1]{PR} carries over to
this more general situation. In particular we have the relation
\begin{equation}
\label{rela}
\wt{Q}_{a,b}(X)=\wt{Q}_{a,b}(X')+
x_1(\wt{Q}_{a-1,b}(X')+\wt{Q}_{a,b-1}(X'))+
x_1^2\wt{Q}_{a-1,b-1}(X'),
\end{equation}
which is true for any $a\gequ b$. The rest of the argument is the 
same.
\end{proof}

\medskip
\noindent
{\bf Remark.} Proposition \ref{ex4.1} fails when $\l$ is a composition;
in fact (\ref{rela}) shows that
it is false even in the case of two part compositions.

\subsection{Symplectic Schubert polynomials}
We denote the elements of the Weyl group $W_n$ for the root system $C_n$ 
as barred permutations, following \cite[\S 1.1]{KT}.
$W_n$ is generated by the elements $s_0,\ldots,s_{n-1}$: for $i>0$,
$s_i$ is the transposition interchanging $i$ and $i+1$, and we define
$s_0$ by
\[
(u_1,u_2,\ldots,u_n)s_0=(\ov{u}_1,u_2,\ldots,u_n).
\]
For each $\l\in \D_n$ of length $\ell$, with $k=n-\ell$ and
$\l'=\rho_n\ssm\l$, the barred permutation 
\[
w_{\l}=(\ov{\l}_1,\ldots,\ov{\l}_{\ell},\l'_k,\ldots,\l'_1)
\]
is the {\em maximal Grassmannian element} of $W_n$ corresponding
to $\l$. We also define the elements
\begin{equation}
\label{wdefs}
w'_{\l}=w_{\l}s_0 \quad \mathrm{and} \quad
w_{\l}''=w_{\l}s_0s_1s_0.
\end{equation}

The group $W_n$ acts on the polynomial ring $\Z[X]$; in particular
the transposition $s_1$ interchanges $x_1$ and $x_2$, 
while $s_0$ replaces $x_1$ by $-x_1$
(all other variables remain fixed). We require the
divided difference
operators $\partial_0, \,\partial'_1: \Z[X]\ra \Z[X]$, defined by
\[
\partial_0(f)=(f-s_0f)/(2x_1)\quad \mathrm{and}
\quad \partial_1'(f)=(f-s_1f)/(x_2-x_1).
\]
Following \cite{PR} and \cite{LP1}, for each $w\in W_n$ there is a 
{\em symplectic Schubert polynomial} $\CS_w(X)\in\Z[X]$, which 
represents the Schubert class associated to $w$ in the cohomology
ring of $Sp_{2n}/B$. According to \cite[Thm.\ A.6]{LP1}, for each
$\l\in\D_n$, the polynomial $\CS_{\l}(X):=\CS_{w_{\l}}(X)$ is equal to
the $\wt{Q}$-polynomial $\wt{Q}_{\l}(X)$. Let us define 
\[
\CS_{\l}'(X)=\CS_{w_{\l}'}(X)
\quad \mathrm{and}\quad
\CS_{\l}''(X)=\CS_{w_{\l}''}(X);
\]
it follows from (\ref{wdefs}) that
\begin{equation}
\label{Cdefs}
\CS_{\l}'(X)=\partial_0(\wt{Q}_{\l}(X))
\quad \mathrm{and}\quad
\CS_{\l}''(X)=\partial_0\partial_1'\partial_0(\wt{Q}_{\l}(X))
\end{equation}
(valid for $\ell(\l)\gequ 1$ and $\ell(\l)\gequ 2$, respectively).

\subsection{Pfaffian identities}
\label{pfids}
We next study certain identities for symplectic Schubert polynomials
which are needed in our
proof of the quantum Giambelli formula for $LG(n,2n)$. We 
see from (\ref{extn}) and (\ref{Cdefs}) that for each strict
partition $\l\in\D_n$,
\begin{equation}
\label{key}
\CS'_{\l}(X)=
\sum_{\substack{k=0 \\ k \ \mathrm{odd}}}^{\ell(\l)}
x_1^{k-1}\sum_{\mu}\wt{Q}_{\mu}(X')
\end{equation}
with the inner sum over all partitions $\mu$ such that
$|\l|-|\mu|=k$ and $\l_i-\mu_i\in\{0,1\}$ for each $i$.
For example, we have
\begin{align}
\label{ex2}
\begin{split}
\CS_{a,b}'(X) &= \wt{Q}_{a-1,b}(X')+\wt{Q}_{a,b-1}(X') \medskip \\
  &= \wt{Q}_{a-1}(X')\wt{Q}_b(X')-
\wt{Q}_a(X')\wt{Q}_{b-1}(X')
\end{split}
\end{align}
for all $a$, $b$ with  $a>b\gequ 0$. 
Our aim is to prove the two Theorems that follow.

\begin{thm}
\label{pfid1} 
Fix $\l\in\D_n$ of length $\ell\gequ 3$, and set
$r=2\lfloor(\ell+1)/2\rfloor$. Then 
\begin{equation}
\label{pf1}
\sum_{j=1}^{r-1}(-1)^{j-1}\,\CS'_{\l_j,\l_r}(X)\,
\CS'_{\l\ssm\{\l_j,\l_r\}}(X)=0.
\end{equation}
\end{thm}

\begin{proof}
Our argument can be conveniently expressed using only
the partitions which index the 
polynomials involved; we therefore begin by defining an algebra 
with formal variables which represent these indices.
Let ${\mathcal A}$ be a commutative
$\Z$-algebra generated by symbols 
$(a_1,a_2,\ldots)$, where the entries $a_i$ are 
barred integers. The symbol $(a_1,a_2,\ldots)$
corresponds to
the polynomial $\wt{Q}_{\rho}(X')$, where $\rho$ is the composition
with $\rho_i$ equal to $a_i$, when $a_i$ is unbarred, and $a_i-1$
when $a_i$ is barred (in the proof, this is only used for 
partitions $\rho$). We identify $(a,0)$ with $(a)$.

Let $\mu$ be a barred partition, that is, a partition in which bars
have been added to some of the entries; assume further
that the underlying partition
is strict. For $\ell(\m)\gequ 3$,
we impose the Pfaffian relation
\begin{equation}
\label{pfrl}
(\m)=
\sum_{j=1}^{r-1}(-1)^{j-1}(\m_j,\m_r)
\cdot (\m\ssm\{\m_j,\m_r\})
\end{equation}
(this corresponds to (\ref{pfaff}) for $\n=\m$). Iterating this
relation gives
\begin{equation}
\label{defin}
(\m)=\sum \epsilon(\m,\nu)(\nu_1,\nu_2)\cdots(\nu_{r-1},\nu_r),
\end{equation}
where the sum is over all $(r-1)(r-3)\cdots (1)$ ways to write the
set $\{\m_1,\ldots,\m_r\}$ as a union of pairs
$\{\nu_1,\nu_2\}\cup\cdots\cup\{\nu_{r-1},\nu_r\}$, and where 
$\epsilon(\m,\nu)$ is the sign of the permutation that takes
$(\m_1,\ldots,\m_r)$ into $(\nu_1,\ldots,\nu_r)$; we adopt the
convention that $\nu_{2i-1}\gequ\nu_{2i}$.

We define the
square bracket symbols $[\m]$, where $\m$ is a strict {\em unbarred}
partition, as follows:
$[a]=(\ov{a})$, $[a,b]=(\ov{a},b)+(a,\ov{b})$, and generally
\begin{equation}
\label{pfaffsym}
[\m]=
\sum_{j=1}^{r-1}(-1)^{j-1}[\m_j,\m_r]
\cdot [\m\ssm\{\m_j,\m_r\}]
\end{equation}
Finally we impose the following relation in our algebra:
\begin{equation}
\label{fghi}
[a,b]=(\ov{a})(b)-(a)(\ov{b});
\end{equation}
this corresponds to (\ref{ex2}).

\begin{lemma}
\label{p1}
For each $\m\in\D_n$ with $\ell(\m)\gequ 3$, we have $[\m]=0$
in ${\mathcal A}$.
\end{lemma}
\begin{proof} When $\ell(\m)\in\{3,4\}$ the identity is easy to check 
using (\ref{fghi}).
For $\m$ of larger length it follows from the Laplace-type expansion
for Pfaffians displayed in (\ref{pfaffsym}).
\end{proof}

\medskip
For a partition $\n$ and integer $k$, we define $B(\nu,k)$ to be
the set of all barred 
partitions $\mu$ which are obtained from $\nu$ by
adding bars over $k$ distinct entries. Now,
using equations (\ref{key}) and (\ref{ex2})
we see that (\ref{pf1}) corresponds to the following
identity, which must hold for each odd $k>0$:
\begin{equation}
\label{idneed}
\sum_{j=1}^{r-1}(-1)^{j-1} [\l_j,\l_r]\cdot
\sum_{\mu\in B(\l\ssm\{\l_j,\l_r\},k)}(\mu)=0.
\end{equation}

To prove (\ref{idneed}), we work as follows. Call a pair 
$[\nu_{2i-1},\nu_{2i}]$ in square brackets {\em distinguished};
pairs which are not distinguished are called {\em normal}.
A part of $\l$ is distinguished if it belongs to a distinguished
pair. We say that a pair $(\nu_{2i-1},\nu_{2i})$ is {\em odd}
if exactly one of the parts $\nu_{2i-1}$ and $\nu_{2i}$ is barred.
In the sequel all brackets $\{\ ,\ \}$ will be either
round $(\ ,\ )$  or square $[\ ,\ ]$, so all expressions are
elements of the algebra ${\mathcal A}$. 

We define a {\em matching} of a barred partition
$\m$ to be an expression of the set
$\{\m_1,\m_2,\ldots,\m_r\}$ as a union of pairs 
$\{\n_1,\n_2\}\cup \cdots \cup \{\n_{r-1},\n_r\}$ as above,
together with a specification of each $\{\n_{2i-1},\n_{2i}\}$
as distinguished or normal.
Fix an odd integer $k$ with $0<k<r-1$. The first claim is that
the left hand side $S$ of the sum (\ref{idneed}) satisfies
\begin{equation}
\label{S}
S=\sum_{\m}\sum_{\nu} 
\epsilon(\l,\nu)\{\nu_1,\nu_2\}\cdots\{\nu_{r-1},\nu_r\}
\end{equation}
where (i) the outer sum 
is over $\mu\in B(\l,k)$ with the part equal to
$\l_r$ unbarred, (ii) 
the inner sum is over all matchings $\nu$ of $\mu$
with a unique distinguished pair, which contains $\l_r$
(iii) the sign
$\epsilon(\l,\nu)$ is defined as above, by ignoring the bars.

Let $t(\m,\n)$ be the total number of odd pairs in the summand 
\[
\epsilon(\l,\nu)\{\nu_1,\nu_2\}\cdots\{\nu_{r-1},\nu_r\}
\]
of $S$; note that $t(\m,\n)>0$. Observe, for each $t$,
that the set of summands with $t(\m,\n)=t$ can be partitioned
into subsets of size $2^t$: elements 
in a single subset $J$ differ only
by moving some of the bars within the odd pairs.
The sum of the elements in $J$ is equal to a single expression
$\epsilon(\l,\nu)\{\nu_1,\nu_2\}\cdots\{\nu_{r-1},\nu_r\}$
with $t$ fewer barred parts, $t+1$ distinguished pairs and no odd pairs.

Working in this way we may eliminate all odd pairs from the sum $S$.
Hence $S=S'$, where
\begin{equation}
\label{S2}
S'=\sum_{\mu}^{\sim}\sum_{\nu}^{\sim}
 \epsilon(\l,\nu)\{\nu_1,\nu_2\}\cdots\{\nu_{r-1},\nu_r\};
\end{equation}
here (i) the outer sum is over all
barred partitions $\mu$ with parts obtained 
by adding an even number, less than $k$, of bars to the parts
of $\l$ 
(ii) the inner sum is over all matchings $\n$ of $\mu$ such that
$\l_r$ is distinguished and the number of barred parts plus
the number of distinguished pairs equals $k+1$,
(iii) the sign $\epsilon(\l,\nu)$ is defined by ignoring the bars, 
as before.

Observe next that $S'$ splits into subsums $S''$; each $S''$ is 
the sum of all summands in $S'$ with a given set of distinguished
{\em parts} and normal {\em pairs}. It is clear that, up to a sign depending
on $S''$, the remainder
after factoring out the normal pairs from $S''$ is a Pfaffian sum like
(\ref{defin}) which involves only distinguished pairs. Now Lemma \ref{p1}
implies that $S''=0$; we deduce that $S'$, and hence $S$, vanishes.
\end{proof}


\medskip

\begin{thm}
\label{pfid2}
For every $\l\in\D_n$ of even length $\ell\gequ 4$ we have
\begin{equation}
\label{pf2}
\sum_{j=1}^{\ell-1} (-1)^{j-1}\,\CS''_{\l_j,\l_{\ell}}(X)\,
\CS''_{\l\ssm\{\l_j,\l_{\ell}\}}(X)=0.
\end{equation}
\end{thm}

\begin{proof}
Using (\ref{rela}) one computes that for $a>b>0$,
\begin{align*}
\partial_0 \partial_1' \partial_0 ( \wt{Q}_{a,b} (X)) &=
\wt{Q}_{a-2,b-1}(X'')+\wt{Q}_{a-1,b-2}(X'') \medskip \\
 &= \wt{Q}_{a-2}(X'')\wt{Q}_{b-1}(X'') - 
     \wt{Q}_{a-1}(X'')\wt{Q}_{b-2}(X''),
\end{align*}
where $X''=(x_3,\ldots,x_n)$. We deduce that 
\[
\mathrm{Pfaffian}[\CS''_{\l_i,\l_j}(X)]_{1\lequ i<j\lequ \ell}=0
\]
for any partition $\l$ of even length $\ell\gequ 4$, by arguing as in 
Lemma \ref{p1}. The proof of the general case of (\ref{pf2})
will require significantly more work.

For each $r$ and $s$ with $r\gequ s\gequ 0$,
let $m_{r,s}(x_1,x_2)$ denote the monomial
symmetric function in $x_1$ and $x_2$. In other words,
$m_{r,s}(x_1,x_2)=x_1^rx_2^s+x_1^sx_2^r$ if $r\neq s$ and
$m_{r,s}(x_1,x_2)=x_1^rx_2^r$ if $r=s$. 

\begin{prop}
\label{lem2}
For any strict partition $\l$ of even length $\ell>0$ we have
\begin{equation}
\label{l2}
\CS''_{\l}(X) =
\sum_{\substack{0\lequ s\lequ r < \ell \\ r,s \ \mathrm{even}}} 
m_{r,s}(x_1,x_2) 
\sum_{\substack{a+2b=r+s+3 \\ a,b\gequ 0}}
\binom{a-1}{s+1-b}\sum_{\nu\in C(\l,a,b)} \wt{Q}_{\nu}(X''),
\end{equation}
where $C(\l,a,b)$ is defined as
the set of compositions $\nu$ with $\l_i-\nu_i\in\{0,1,2\}$ for 
all $i$ and $\l_i-\nu_i=1$ {\em (}resp.\ $\l_i-\nu_i=2${\em )}
for exactly $a$ {\em (}resp.\ $b${\em )} values of $i$.
\end{prop}
\begin{proof} 
Use (\ref{Cdefs}) and (\ref{key}) to compute
\begin{equation}
\label{cequ}
\CS''_{\l}(X)=\partial_0 \partial_1' 
\Bigl(\sum_{\substack{j=0 \\ j \ \mathrm{even}}}^{\ell-2}
x_1^j\sum_{\mu\in P(\l,j+1)}\wt{Q}_{\mu}(X')\Bigr),
\end{equation}
where $P(\l,j+1)$ is the set of compositions 
$\mu$ with $|\l|-|\mu|=j+1$ and
$\l_i-\mu_i\in\{0,1\}$ for each $i$. Note also that for 
any $j,k\gequ 0$ we have
\begin{equation}
\label{d1}
\partial_1 (x_1^jx_2^k)=\mathrm{sgn}(j-k)
\sum_{\substack{c+d=j+k-1 \\ c,d\gequ\min\{j,k\}}}x_1^cx_2^d.
\end{equation}

For each strict partition
$\l$ and $j,k\gequ 0$, let 
$P(\l,j,k)$ be the set of pairs of compositions $(\mu,\nu)$ such that
$\mu\in P(\l,j)$ and $\nu\in P(\mu,k)$.
Now apply (\ref{extn}) and (\ref{d1}) in (\ref{cequ}) to get
\begin{align*}
\CS''_{\l}(X)&=\sum_{j \ \mathrm{even}}
\partial_0 \partial_1'(x_1^jx_2^k)
\sum_{(\m,\n) \in P(\l,j+1,k)}\wt{Q}_{\nu}(X'')\\
&=\sum_{j \ \mathrm{even}}
\mathrm{sgn}(k-j)\sum_{\substack{ c+d=j+k-1 \\ c,d\gequ \min\{j,k\} \\
c\ \mathrm{odd}}}x_1^{c-1}x_2^d
\sum_{(\m,\n) \in P(\l,j+1,k)}\wt{Q}_{\nu}(X'').
\end{align*}
It is useful to
observe that $\CS''_{\l}(X)$ is a polynomial in $R[x_1^2,x_2^2]^{S_2}$,
where $R=\Z[x_3,\ldots,x_n]^{S_{n-2}}$. This follows because 
\[
\partial_0(\CS''_{\l}(X))=\partial_1(\CS''_{\l}(X))=0
\]
(for the latter, use the fact that 
$\partial_1\partial_0\partial_1\partial_0=
\partial_0\partial_1\partial_0\partial_1$ and $\partial_1(\wt{Q}_{\l}(X))=0$).
Using this to cancel all irrelevant terms in the previous equality gives
\begin{equation}
\label{finally}
\CS''_{\l}(X)=\sum_{\substack{0\lequ s\lequ r < \ell \\ r,s \ \mathrm{even}}} 
m_{r,s}(x_1,x_2)
\sum_{\substack{j,k \ \mathrm{even}
\\ j+k=r+s+2 \\\min\{j,k\}\lequ r,s}}
\sum_{(\m,\n) \in P(\l,j+1,k)}
\mathrm{sgn}(k-j)\,\wt{Q}_{\nu}(X'').
\end{equation}

Consider the terms occuring in the inner sum $\sum\sum$ of (\ref{finally})
which involve a fixed $\nu\in C(\l,a,b)$, for some $a,b\gequ 0$
such that $a+2b=r+s+3$. Observe that the coefficient 
of any such term depends only on $(r,s,a,b)$ and
can be expressed as
\begin{equation}
\label{binsum}
\sum_{\substack{j,k \ \mathrm{even}
\\ j+k=r+s+2 \\\min\{j,k\}\lequ r,s}}\mathrm{sgn}(k-j)
\binom{a}{k-b}.
\end{equation}
Now (\ref{binsum}) can be 
rearranged as an alternating sum of binomial coefficients,
which in turn simplifies to the single binomial coefficient
displayed in (\ref{l2}).
\end{proof}

\medskip

For any $r$ and $s$, both even, the coefficient of $m_{r,s}(x_1,x_2)$
in the left-hand side of (\ref{pf2}), expanded as a polynomial in
$R[x_1^2,x_2^2]^{S_2}$, is
\begin{equation}
\label{zerosum}
\sum_{j=1}^{\ell-1} (-1)^{j-1}
\CS''_{\l_j,\l_{\ell}}(X)
\sum_{\substack{a+2b=r+s+3 \\ a,b\gequ 0}} 
\binom{a-1}{s+1-b}
\sum_{\nu\in C(\l\ssm\{\l_j,\l_{\ell}\},a,b)} 
\wt{Q}_{\nu}(X'').
\end{equation}
Hence, to prove Theorem \ref{pfid2} it suffices to show that the
expression (\ref{zerosum}) vanishes for any strict partition $\l$
of even length $\ell\gequ 4$ and any $r,s\gequ 0$, both even.

To accomplish this, we form
an algebra ${\mathcal B}$ as in the proof of Theorem
\ref{pfid1}, except that the symbols $(a_1,a_2,\ldots)$ are now
sequences of integers which can each have up to two bars. 
This time the symbol $(a_1,a_2,\ldots)$
corresponds to
the polynomial $\wt{Q}_{\n}(X'')$, where $\n$ is the composition
with $\n_i$ equal to the integer $a_i$ minus the number of bars
over $a_i$. 

As before, we impose the Pfaffian relations (\ref{pfrl}) in
${\mathcal B}$. Define also $[a,b]=(\ov{a},b)+(a,\ov{b})$,
where $a$ and $b$ are integers, each with up to one bar. More
generally, define $[\m]$ for any barred partition $\m$, as in
(\ref{pfaffsym}). We also impose the relations
\[
[a,b]=(\ov{a})(b)-(a)(\ov{b}) 
\]
for integers $a$, $b$, with up to one bar each.
As in the proof of Theorem \ref{pfid1}, we have
\begin{lemma}
\label{flem}
Let $\m$ be a strict partition of even length at least $4$,
such that either $\m$ is unbarred or every part of $\mu$ has
a bar. Then we have $[\m]=0$ in ${\mathcal B}$.
\end{lemma}


The vanishing of (\ref{zerosum}) corresponds to the following identity
in ${\mathcal B}$:
\begin{equation}
\label{Bzsum}
\sum_{\substack{a+2b=r+s+3 \\ a,b\gequ 0}} 
\binom{a-1}{s+1-b}
\sum_{j=1}^{\ell-1} (-1)^{j-1}
[\ov{\l}_j,\ov{\l}_{\ell}]
\sum_{\nu\in C(\l\ssm\{\l_j,\l_{\ell}\},a,b)} 
(\nu) = 0.
\end{equation}
For each $a$ and $b$, denote by $W_{a,b}$ the contents of the inner
two sums in (\ref{Bzsum}). 
By applying the Pfaffian expansion to the terms $(\nu)$ in $W_{a,b}$,
we can write 
\[
W_{a,b}=\sum_{\rho} 
\epsilon(\l,\rho)\{\rho_1,\rho_2\}\cdots\{\rho_{\ell-1},\rho_{\ell}\}.
\]
The summands in $W_{a,b}$ are of three types: first, those that contain
a pair $(\rho_{2i-1},\rho_{2i})$ with a total of three bars; second, those
not of the first type which contain at least two pairs, each 
with exactly one bar; third, the remaining summands. Let $R_{a,b}$,
$S_{a,b}$ and $T_{a,b}$ denote the sum of the summands of, respectively,
the first, second and third types.

Arguing as in the proof of Theorem \ref{pfid1} and using Lemma \ref{flem},
we see that $R_{a,b}=S_{a,b}=0$, for all $a$, $b$. However, it is not true
that $T_{a,b}$ is always zero.
Observe that $T_{a,b}$ splits as a sum $\sum_gT^g_{a,b}$ indexed by
the number $g$ of pairs in each summand which contain $4$ bars. 
We will show that
\begin{equation}
\label{needshow}
\sum_{\substack{a+2b=r+s+3 \\ a,b\gequ 0}} 
\binom{a-1}{s+1-b} T^g_{a,b}=0
\end{equation}
for each $g$. In the following we give the argument when $g=0$, and
afterwards we describe the small modifications for the proof of the
general case.

Assume that $r\gequ s$ and introduce 
\[
u=(r-s)/2 \ \ \ \ \ \mathrm{and} \ \ \ \ \ v=(r+s)/2.
\]
We recursively define a sequence of coefficients: Put $e_u=1$ and
\[
e_m=\binom{2m}{m-u}-\frac{2m}{v+2-m}\,e_{m-1}
\]
for $m=u+1,\ldots,v+1$. 

\begin{lemma}
\label{lastzero}
We have $e_{v+1}=0$.
\end{lemma}
\begin{proof}
For any integer $p\gequ 0$ we have the combinatorial identity
\begin{equation}
\label{gendawson}
\sum_k(-1)^k2^{-k}\binom{p}{k}\binom{2k}{k-q}=
\begin{cases}
(-1)^q2^{-p}\binom{p}{(p+q)/2}&\text{if $p+q$ is even},\\
0&\text{if $p+q$ is odd}.
\end{cases}
\end{equation}
This is proved, e.g., by showing the left-hand side $L(p,q)$ satisfies
the recursion $2L(p,q)+L(p-1,q-1)+L(p-1,q+1)=0$ and inducting on $p$.
When $q=0$, (\ref{gendawson})
reduces to an identity attributed to Dawson in \cite[p.\ 71]{R}.
The assertion $e_{v+1}=0$ follows from the case
$(p,q)=(v+1,u)$ of (\ref{gendawson}).
\end{proof}

\medskip
We deduce from Lemma \ref{lastzero} that
\begin{multline*}
\sum_{\substack{a+2b=r+s+3 \\ a,b\gequ 0}} 
\binom{a-1}{s+1-b} T^0_{a,b} =
\sum_{m=u}^{v+1}\binom{2m}{m-u}T^0_{2m+1,v+1-m} \\
=
\sum_{m=u+1}^{v+1}e_{m-1}
\Bigl(T^0_{2m-1,v+2-m}+\frac{2m}{v+2-m}\,T^0_{2m+1,v+1-m}\Bigr).
\end{multline*}

\medskip

\begin{lemma}
\label{lem4}
We have 
\begin{equation}
\label{tzeroident}
(v+2-m)\,T^0_{2m-1,v+2-m}+2m\,T^0_{2m+1,v+1-m}=0
\end{equation}
for every integer $m$ with $u+1\lequ m\lequ v+1$.
\end{lemma}
\begin{proof}
We have that
\[
T^0_{a,b}=\sum_{\m\in C(\l,a+2,b)}\sum_{\nu} 
\epsilon(\l,\nu)\{\nu_1,\nu_2\}\cdots\{\nu_{\ell-1},\nu_{\ell}\},
\]
where the inner sum is over all matchings $\nu$ of $\mu$
with (i) a unique distinguished pair, which contains $\l_{\ell}$ and has
bars on both entries, (ii) a unique pair which contains only one bar, and
(iii) no pairs which contain $3$ or $4$ bars. 

Now,
\begin{equation}
\label{feq}
(v+2-m)\,T^0_{2m-1,v+2-m}=
\sum_{\m\in C(\l,2m+1,v+2-m)}\sum_{(\nu,i)} 
\epsilon(\l,\nu)\{\nu_1,\nu_2\}\cdots\{\nu_{\ell-1},\nu_{\ell}\},
\end{equation}
where the inner sum is over pairs $(\n,i)$, where $\n$ is a matching
satisfying (i)--(iii), and $i$ is the index of a double-barred part
of $\n$. Also, 
\begin{equation}
\label{seq}
2m\,T^0_{2m+1,v+1-m}=
\sum_{\m\in C(\l,2m+3,v+1-m)}\sum_{(\rho,j)} 
\epsilon(\l,\rho)\{\rho_1,\rho_2\}\cdots\{\rho_{\ell-1},\rho_{\ell}\},
\end{equation}
where the inner sum is over pairs $(\rho,j)$, where $\rho$ is a matching
satisfying (i)--(iii),
and where $j$ is the index of a single-barred part of $\rho$,
which is not in the distinguished pair, nor in the 
unique pair which contains only one bar. Observe that there is a bijection
between the pairs $(\n,i)$ which appear in (\ref{feq}) and pairs
$(\rho,j)$ in (\ref{seq}).

Adding (\ref{feq}) and (\ref{seq}) and using the bijection between the
indexing sets and the identity
\[
(\ov{\ov c}, d) + 2(\ov c, \ov d) + (c,\ov{\ov d}) =
[\ov c, d] + [c,\ov d],
\]
we may express the left-hand side of (\ref{tzeroident}) as a sum in
which every summand has two distinguished pairs.
By invoking the identity
\begin{multline*}
[\ov a, \ov b]([\ov c, d]+[c,\ov d])
-[\ov a, \ov c]([\ov b, d]+[b,\ov d])
+[\ov a,\ov b]([\ov b, c]+[b,\ov c])    \\
=
-(\ov{\ov a})(\ov{\ov b})[c,d]
+(\ov{\ov a})(\ov{\ov c})[b,d]
-(\ov{\ov a})(\ov{\ov d})[b,c]
\end{multline*}
(with $a=\l_{\ell}$, always),
we reduce ourselves to a situation where we may use
Lemma \ref{flem} and conclude that (\ref{tzeroident}) holds.
\end{proof}

\medskip
It follows from Lemma \ref{lem4} and the previous discussion that
(\ref{needshow}) holds for $g=0$.
When $g$ is positive, the above argument goes through with the triple
$(b,r,s)$ replaced with $(b-2g, r-2g, s-2g)$; note that the binomial
coefficients in (\ref{needshow}) remain unchanged under these
substitutions. This finishes the proof of Theorem \ref{pfid2}.
\end{proof}

\section{Isotropic Grassmannians}
\label{isograss}

\subsection{Schubert subvarieties and incidence loci}
\label{schuvar}

The principal object of study is the Lagrangian Grassmannian
$LG(n,2n)$ which parametrizes $n$-dimensional subspaces of a fixed
$2n$-dimensional complex vector space $V$, isotropic for a fixed
nondegenerate skew-symmetric bilinear form on $V$.
When $n$ is fixed, we write $LG$ for $LG(n,2n)$.
We have $\dim_{\C}LG=n(n+1)/2$. The identities in cohomology that we 
establish
in \S \ref{schuvar} and \S \ref{pfaffident} remain valid if we work
over an arbitrary base field, using Chow rings in place of cohomology.

Let $F_{\bullet}$ be a fixed complete isotropic flag of subspaces of $V$.
The Schubert varieties $\X_\l\subset LG$ --- 
defined relative to $F_{\bullet}$ --- are indexed by
partitions $\l\in \D_n$. 
One has
$\dim(\Sigma\cap F_k^\perp)=\dim(\Sigma\cap F_k)+n-k$
for any Lagrangian subspace $\Sigma$, so there are two ways to
write the conditions which define a Schubert variety in $LG$:
\begin{align}
\X_\l &= \{\,\Sigma\in LG\,|\,
\rk(\Sigma\to V/F_{n+1-\l_i})\lequ n-i,\,\,i=1,\ldots,\ell(\l)\,\}
\label{schubalt} \\
&= \{\,\Sigma\in LG\,|\,
\rk(\Sigma\to V/F^\perp_{n+1-\l_i})
\lequ n+1-i-\l_i,\,\,i=1,\ldots,\ell(\l)\,\}.
\label{schubert}
\end{align}
We call particular attention to the alternative formulation (\ref{schubert}),
as these are the relevant rank conditions for situations where the morphism
$\Sigma\to V$ is allowed to degenerate.

Set $\s_{\l}=[\X_{\l}]$ in $H^*(LG,\Z)$.
The classical Giambelli formula (\ref{lgqgiam2})
for $LG$ is equivalent to the following identity in $H^*(LG,\Z)$:
\begin{equation}
\label{classicalgiambelli}
\s_\l = \sum_{j=1}^{r-1} (-1)^{j-1} \s_{\l_j,\l_r}\cdot
\s_{\l\smallsetminus\{\l_j,\l_r\}},
\end{equation}
for $r=2\lfloor (\ell(\l)+1)/2 \rfloor$.
For $\m\in \D_n$, let us denote $\mu'=\rho_n\ssm\mu$, the dual partition.
Then, we recall, the Poincar\'e duality pairing on $LG$ satisfies
\[
\int_{LG} \s_{\l}\,\s_{\m} = \delta_{\l\m'}.
\]

Later on,
in place of $\Sigma\to V$, we will have
a morphism of vector bundles $\E\to \cO_T\otimes V$ over some base $T$,
with $\E$ a rank $n$ vector bundle and with the morphism generically of
full rank, but with loci where the rank drops.
So we need to study the Grassmannians $IG(k,2n)$ of isotropic $k$-dimensional
subspaces of $V$, for various $k<n$, notably for $k=n-1$ and $k=n-2$.
We have $IG(n,2n)=LG(n,2n)$.

Observe that any $Sp_{2n}$-translate of the Schubert variety 
$\X_{n,n-1,\ldots,k+1}$ in $LG$ is of the form
$\{\,\Sigma\in LG \ |\ \Sigma\supset A\,\}$
for a unique $A\in IG(n-k,2n)$.
Any such translate can be identified with $LG(k,2k)$, and moreover,
for any $\l$, meets $\X_{\l}$ in a
Schubert subvariety of $LG(k,2k)$:

\begin{prop}
\label{meetsmallergrass}
Let $A$ be an isotropic subspace of $V$ of dimension $n-k$, and 
let $Y\subset LG(n,2n)$ be the subvariety of
Lagrangian subspaces of $V$ which contain $A$.
Then $\X_\lambda\cap Y$ is a Schubert variety in $Y\simeq LG(k,2k)$ for
any $\l\in\D_n$. Moreover, if $\ell(\l)<k$ then the intersection,
if nonempty, has positive dimension.
\end{prop}

\begin{proof}
Define the isotropic flag $\wt{F}_\bullet$ of subspaces of $A^\perp/A$
by $\wt{F}_i=((F_i+A)\cap A^\perp)/A$.
For any Lagrangian $\Sigma\subset V$ containing $A$, we have
$\dim(\Sigma\cap F_i)=\dim((\Sigma/A)\cap \wt{F}_i)+\dim(A\cap F_i)$.
So $\X_\l\cap Y$ is defined by the attitude of $\Sigma/A$ with respect
to $\wt{F}_\bullet$, and hence is a Schubert variety (if nonempty). 
For the intersection to be a point, we would require at least $k$
rank conditions as in (\ref{schubalt}), that is, $\ell(\l)\gequ k$.
\end{proof}

\medskip
The space $IG(n-1,2n)$ is the parameter space of lines on $LG$, and
the variety of lines incident to $\X_\l$ (for nonempty $\l$)
is the Schubert variety
\begin{equation}
\label{schubertprime}
\X'_\l = \{\,\Sigma'\in IG(n-1,2n)\,|\,
\rk(\Sigma'\to V/F^\perp_{n+1-\l_i})\lequ n+1-i-\l_i,\,\,{\rm all}\ i\,\}.
\end{equation}
Note the rank conditions are identical to those in (\ref{schubert}).
The codimension of $\X'_\l$ is $|\l|-1$.
Analogously, $IG(n-2,2n)$ parametrizes translates of $\X_{n\cdots43}$ 
(isomorphic to the quadric threefold $LG(2,4)$) on $LG$.
Now suppose $\ell(\l)\gequ 2$;
then translates of $\X_{n\cdots43}$
on $LG$ incident to $\X_\l$ form the
Schubert variety
\begin{equation}
\label{schubertprimeprime}
\X''_\l = \{\,\Sigma''\in IG(n-2,2n)\,|\,
\rk(\Sigma''\to V/F^\perp_{n+1-\l_i})\lequ n+1-i-\l_i,\,\,{\rm all}\ i\,\},
\end{equation}
which has codimension $|\l|-3$. We remark that 
the subvarieties we have described 
are only some of the Schubert varieties in $IG(n-1,2n)$ and 
$IG(n-2,2n)$.

The modular interpretation of the loci $\X'_\l$ and $\X''_\l$ 
is explained by 
\begin{lemma}
\label{observation}
Let $A$ be an isotropic subspace of $V$ of dimension $n-k$, and
let $F_n$ be a fixed Lagrangian subspace of $V$.
Then there exists a unique Lagrangian subspace $\Sigma$ which contains
$A$ and satisfies
$\dim(\Sigma\cap F_n) = \dim (A\cap F_n)+k$. In fact, we have
$\Sigma=\Span(A,A^{\perp}\cap F_n)$.
\end{lemma}
\begin{proof}
Set $W=A\cap F_n$; replacing $V$ by $W^\perp/W$ and
$A$ and $F_n$ by their respective images in $W^\perp/W$,
we are reduced to the case $A\cap F_n=0$.
Now for dimension reasons, we have
$\dim (A^{\perp}\cap F_n)\gequ k$,
so isotropicity of $\Span(A, A^{\perp}\cap F_n)$
forces $\dim (A^{\perp}\cap F_n)=k$,
and the assertion is clear.
\end{proof}

%

\subsection{Pfaffian identities on isotropic Grassmannians}
\label{pfaffident}
Let $F=F_{Sp}(V)$ denote the variety of complete isotropic flags
in $V$. There are natural projection maps from $F$ to
the Grassmannians $IG(n-1,2n)$ and $IG(n-2,2n)$, inducing
injective pullback morphisms on cohomology. Referring to 
\cite[\S 2.4]{KT}, one sees that the Schubert class
$[\X'_\l]$ (resp.\ $[\X''_\l]$) in $H^*(IG(n-1,2n))$ (resp.\ 
$H^*(IG(n-2,2n))$ pulls back to the class $\CS'_{\l}(X)$
(resp.\ $\CS''_{\l}(X)$) in $H^*(F)$, for each $\l\in\D_n$. 
Here $X=(x_1,\ldots,x_n)$ is the vector of Chern roots of the dual
to the tautological rank $n$ vector bundle over $F$, ordered as
in \cite[Sect.\ 2]{KT}. Theorems \ref{pfid1} and \ref{pfid2}
now give

\begin{cor}
\label{forcor} 
{\em a)} For every $\l\in\D_n$ of length $\ell\gequ 3$ 
and $r=2\lfloor(\ell+1)/2\rfloor$ we have
\begin{equation}
\label{idone}
\sum_{j=1}^{r-1}(-1)^{j-1}\,[\X'_{\l_j,\l_r}]\,
[\X'_{\l\ssm\{\l_j,\l_r\}}]=0
\end{equation}
in $H^*(IG(n-1,2n),\Z)$.

\medskip
\noin
{\em b)} For every $\l\in\D_n$ of even length $\ell\gequ 4$ we have
\begin{equation}
\label{idtwo}
\sum_{j=1}^{\ell-1} (-1)^{j-1}\,[\X''_{\l_j,\l_{\ell}}]\,
[\X''_{\l\ssm\{\l_j,\l_{\ell}\}}]=0
\end{equation}
in $H^*(IG(n-2,2n),\Z)$.
\end{cor}

\subsection{Counting lines on $LG(n,2n)$}
\label{countinglines}
In this subsection, we exploit a correspondence between lines on
$LG(n,2n)$ and points on $LG(n+1,2n+2)$ to derive the following formula
for three-point degree-$1$ Gromov--Witten invariants on $LG=LG(n,2n)$:
\begin{prop}
\label{linenumbers}
For $\l$, $\m$, $\n\in \D_n$ we have 
\[
\langle \sigma_\l,\sigma_\m,\sigma_\n\rangle_1=
\frac{1}{2}\int_{LG(n+1,2n+2)}[\X^+_\l]\cdot[\X^+_\m]\cdot[\X^+_\n],
\]
where $\X_\l^+$, $\X^+_\m$, $\X^+_\n$ 
 denote Schubert varieties in $LG(n+1,2n+2)$.
\end{prop}

\begin{proof}
Let $V$ be as in \S \ref{schuvar}, let $H$ be a $2$-dimensional
symplectic vector space and let $V^+$ be the orthogonal direct
sum of $V$ and $H$.
Consider the correspondence between $LG(n+1,2n+2)$ and
$IG(n-1,2n)$ consisting of pairs $(\Sigma^+,\Sigma')$
with $\Sigma^+$ a Lagrangian subspace of $V^+$ and $\Sigma'$
an isotropic $(n-1)$-dimensional subspace of $V$, given by the
condition $\Sigma'\subset\Sigma^+$.
This is the correspondence induced by the rational map which sends
$[\Sigma^+]$ to $[\Sigma^+\cap V]$,
defined for Lagrangian $\Sigma^+\subset V^+$
with $\dim(\Sigma^+\cap V)=n-1$. Choose a nonzero $h\in H$ and
consider the isotropic flags $F_{\bullet}$
in $V$ and $F^+_{\bullet}$ in $V^+$, with $F_i^+=F_{i-1}\oplus\langle h
\rangle$. Then, for any $\l\in \D_n$, we have that to any $\Sigma^+$
with $\Sigma^+\in \X_\l^+\subset LG(n+1,2n+2)$,
there corresponds $\Sigma'$ with $\Sigma'\in \X'_\l$.

The degree of $q$ in $QH^*(LG(n,2n))$ equals
\[
\int_{LG}c_1(T_{LG})\cdot \s_{1'} = n+1,
\]
so for degree reasons, the Gromov--Witten invariant
$\langle \sigma_\l,\sigma_\m,\sigma_\n\rangle_1$
is potentially nonzero (and hence counts lines on $LG$) when
\[
|\l|+|\m|+|\n|=(n+1)(n/2+1).
\]
Now choose general translates of isotropic flags of subspaces of $V$
such that $\X'_\l$, $\X'_\m$, and $\X'_\n$ meet transversely, and such that
every point in the intersection corresponds to a line
in $LG(n,2n)$ which

(i) is incident to $\X_\l$, $\X_\m$, and $\X_\n$ at three distinct points;

(ii) has each incidence point in the corresponding Schubert cell
inside the

\ \ \ \ \ Schubert variety (i.e., all dimension estimates are sharp);

(iii) (for $\ell=\ell(\l)<n$) the line is not incident to
$\X_{\rho_{\ell+1}}$, and similarly for $\m$ and

\ \ \ \ \ $\nu$ (using the respective defining flags).

\noindent
We remark that (ii) and (iii) can be combined into the single statement,
that the points in $\X'_\l\cap \X'_\m\cap \X'_\n$ all lie in the intersections
of the corresponding Schubert cells in $IG(n-1,2n)$.

For each of $\l$, $\m$, $\n$, we extend the corresponding flag of
isotropic spaces in $V$ to subspaces of $V^+$ by adjoining, in each case,
a generally chosen element of $H$.
Now it is evident that every point in $\X^+_\l\cap \X^+_\m\cap \X^+_\n$
(intersection on $LG(n+1,2n+2)$) corresponds to a point in
$\X'_\l\cap \X'_\m\cap \X'_\n$.
Analysis using Lemma \ref{observation}
shows that to every point in
$\X'_\l\cap \X'_\m\cap \X'_\n$ there correspond
exactly two points in $\X^+_\l\cap \X^+_\m\cap \X^+_\n$, each a point of
transverse intersection.
This uses the fact that on $LG(2,4)$ (case $n=1$),
three general $SL_2\times SL_2$
translates of $\X_1^+$ meet transversely at two points, both away from
the locus of $\Sigma^+$ meeting $V$ nontrivially.
\end{proof}

\begin{cor} If $i$ and $j$ are such that $i+j\gequ n+1$, then
\[
\sigma_i \sigma_j = 2\sum_{k=1}^{n-i} \sigma_{i+k,j-k}
+ \sigma_{i+j-n-1}\, q
\]
in $QH^*(LG(n,2n))$.
\end{cor}

{}From the Corollary we immediately deduce that the quantum relation
(\ref{lgqrel}) and the two-condition quantum Giambelli formula
(\ref{lgqgiam1}) are valid in $QH^*(LG)$.
By \cite{ST}, now, we obtain a presentation of $QH^*(LG)$ as a 
quotient of the polynomial ring $\Z[\s_1,\ldots,\s_n,q]$ modulo the
relations (\ref{lgqrel}) (see also \cite[Sect.\ 10]{FPa}).
The proof of the more difficult relation (\ref{lgqgiam2}) occupies
Sections \ref{lagrquot} and \ref{itlqd}.

\section{Lagrangian Quot schemes}
\label{lagrquot}

\subsection{Overview}
\label{overview}

Here is a summary of our analysis in the next two sections.
The goal is to establish an identity in the Chow group of the Lagrangian Quot
schemes, from which identity (\ref{lgqgiam2}) in $QH^*(LG(n,2n))$ follows.
We make use of type $C$ degeneracy loci for isotropic morphisms
of vector bundles \cite{KT} to define classes
$[W_\l(p)]_k$ ($p\in \bP^1$)
of the appropriate dimension $k:=(n+1)(n/2+d)-|\l|$ in the
Chow group of the
Quot scheme $LQ_d$, which compactifies the space of 
degree-$d$ maps
$\bP^1\to LG$.
Choose a point
$p'$ on $\bP^1$ distinct from $p$, and use $\wt{W}$ to
denote a degeneracy
locus defined with respect to a general translate of the fixed
isotropic flag of vector spaces. 
We produce a Pfaffian formula analogous to (\ref{classicalgiambelli}):
\begin{equation}
\label{Wequ}
[W_{\l}(p)]_k=\sum_{j=1}^{r-1}(-1)^{j-1}[W_{\l_j,\l_r}(p)\cap
\wt{W}_{\l\ssm\{\l_j,\l_r\}}(p')]_k,
\end{equation}
for any $\l\in\D_n$ with $\ell(\l)\gequ 3$ and
$r=2\lfloor(\ell(\l)+1)/2\rfloor$.

In fact, we need the cycles in (\ref{Wequ})
to remain rationally equivalent under
further intersection with some (general translate of) $W_{\m}(p'')$
(with varying $\m$, where $p''$ is a third distinct point on $\bP^1$).
The way around the difficulty of a lack of
an intersection product on the (possibly singular)
$LQ_d$ is to work instead
on a modification $LQ_d(p'')$, on which the
evaluation-at-$p''$ map is globally defined, and employ a
refined intersection operation from $LG$.
Then, the rational equivalence implies that the (zero!)\ number of
degree-$d$ maps $\bP^1\to LG$ which send $p$ to $\X_\l$ and $p''$ to
a general translate of $\X_\m$ --- with no constraint on any third point --- 
is equal to the alternating sum of $3$-point degree-$d$ structure
constants $\langle \s_{\l_j,\l_r}, \s_{\l\ssm\{\l_j,\l_r\}}, \s_{\m}\rangle_d$,
and this is
exactly what we require to establish quantum Giambelli
(Theorem \ref{qgthm}).

The required rational equivalence
(which is similar to (\ref{Wequ}), but involves cycles on $LQ_d(p'')$)
is a sum of rational equivalences coming from three sources:
(i) Pfaffian formulas on $LG$ (\ref{classicalgiambelli});
(ii) supplementary rational equivalences on the isotropic Grassmannians of
$(n-1)$-dimensional and $(n-2)$-dimensional isotropic subspaces
of the ambient vector space (\ref{idone}), (\ref{idtwo});
(iii) rational equivalences $\{p\}\sim\{p'\}$ on $\bP^1$.
The most interesting of these are (ii).
Indeed, at some stage in the analysis,
we consider intersections
$W_{\l_j,\l_r}(p)\cap \wt{W}_{\l\ssm\{\l_j,\l_r\}}(p)$,
where the same point $p$ appears twice.
Such an intersection has several components of the correct dimension.
Of course, degree-$d$ morphisms $\bP^1\to LG$ sending $p$ into
$\X_{\l_j,\l_r}\cap \wt{\X}_{\l\ssm\{\l_j,\l_r\}}$ make up one such component.
But also there are components of the same dimension, entirely
contained in the boundary of the Quot scheme.
The formulas in Corollary \ref{forcor} give us precisely the cancellation
of these extra contributions.

\subsection{Definition of $LQ_d$}
\label{lagrquotdefn}
Fix a complex vector space $V$ and set $N=\dim V$.
Let $m$ and $n$ be positive integers, with $m+n=N$.
We recall (from \cite{G})
that the Quot scheme $Q_d$ ($d\gequ 0$)
parametrizing quotient sheaves of $\cO_{\bP^1}\otimes V$
with Hilbert polynomial $nt+n+d$ is a smooth projective variety
which compactifies the space of parametrized degree-$d$ maps from $\bP^1$
to the Grassmannian of $m$-dimensional subspaces of $V$.
On $\bP^1\times Q_d$ there is a universal exact sequence of sheaves
\[
0\lra \E\lra \cO\otimes V\lra \QQ\lra 0
\]
with $\E$ locally free of rank $m$.

We fix $m=n$ and $N=2n$ from now on, and let
the $2n$-dimensional vector space $V$ be endowed with a nondegenerate
skew-symmetric bilinear form.

\begin{defn}
Let $d$ be a nonnegative integer.
The {\em Lagrangian Quot scheme} $LQ_d$ is
the closed subscheme of $Q_d$ which is defined by the vanishing
of the composite
$$\E\lra \cO_{\bP^1}\otimes V 
\stackrel{\alpha}\lra \cO_{\bP^1}\otimes V^*\lra \E^*$$
where $\alpha$ is the isomorphism defined by the symplectic form on $V$.
\end{defn}

For $n=2$, the scheme $LQ_1$ was introduced and studied in
\cite{KT};
it is nonsingular.
For general $n$ and $d$, however, $LQ_d$ has singularities.
The open subscheme $LM_d\subset LQ_d$, defined as the locus where
$\E\to \cO_{\bP^1}\otimes V$ has everywhere full rank,
is smooth and irreducible, and coincides with the moduli space
$M_{0,3}(LG,d)$ of $3$-pointed (that is, parametrized) maps
$\bP^1\to LG$ of degree $d$ (the irreducibility of all $M_{0,i}(LG,d)$
is known \cite{kimpandharipande} \cite{thomsen}).

\subsection{Degeneracy loci}
\label{degloci}
Degeneracy loci for isotropic morphisms of vector bundles in type $C$
(as well as types $B$ and $D$) were introduced in \cite{KT}.

\begin{defn}
The degeneracy loci $W_\l$ and $W_\l(p)$ ($\l\in \D_n$ and $p\in \bP^1$)
are the following subschemes of $\bP^1\times LQ_d$:
\begin{align*}
W_\l &= \{\,x\in \bP^1\times LQ_d\,|\,
\rk(\E\to \cO\otimes V/F^\perp_{n+1-\l_i})_x\lequ n+1-i-\l_i,\\
&\qquad\qquad\qquad\qquad\qquad\qquad\qquad\qquad\qquad\qquad\qquad
 i=1,\ldots,\ell(\l)\,\}, \\
W_\l(p) &= W_\l\cap (\{p\}\times LQ_d)
\end{align*}
\end{defn}

Define also
$$h(n,d)=(n+1)(n/2+d),$$
the dimension of the Lagrangian Quot scheme $LQ_d$.
The fact that
three-term Gromov--Witten invariants on $LG$ arise as numbers
of points in intersections of degeneracy loci on $LQ_d$ is a consequence
of the
\begin{movlemma}
Let $k$ be a positive integer, and let $p_1$, $\ldots$, $p_k$ be
distinct points on $\bP^1$.
Let $\l^1$, $\ldots$, $\l^k$ be partitions in $\D_n$, and
let us take the degeneracy loci $W_{\l^1}(p_1)$, $\ldots$, $W_{\l^k}(p_k)$
to be defined by isotropic flags of vector spaces in general position.
Let us define
$$Z=W_{\l^1}(p_1)\cap\cdots\cap W_{\l^k}(p_k).$$
Then $Z$ has dimension at most $h(n,d)-\sum_{i=1}^k |\l^i|$.
Moreover, $Z\cap LM_d$ is either empty or generically reduced and
of pure dimension $h(n,d)-\sum_i |\l^i|$;
also, $Z\cap (LQ_d\smallsetminus LM_d)$ has dimension at most
$h(n,d)-\sum_{i=1}^k |\l^i|-1$.
\end{movlemma}

\noindent
We prove the Moving Lemma in \S \ref{structure}, after stating
a theorem which describes the structure of the boundary of $LQ_d$.

\begin{cor}
\label{corgw}
Let $p$, $p'$, $p''\in \bP^1$ be distinct points.
Suppose $\l$, $\m$, $\n\in \D_n$ satisfy
$|\l|+|\m|+|\n|=h(n,d)$.
With degeneracy loci defined with respect to isotropic flags in general
position, the intersection $W_\l(p)\cap W_\m(p')\cap W_\n(p'')$
consists of finitely many reduced points, all contained in $LM_d$, and 
the corresponding Gromov--Witten invariant on $LG$ equals the cardinality of
this set of points:
$$\langle \s_\l, \s_\m, \s_\n \rangle_d = \#\bigl(W_\l(p)\cap W_\m(p')\cap
W_\n(p'')\bigr).$$
\end{cor}

\begin{cor}
\label{corzero}
If $p$ and $p'$ are distinct points of $\bP^1$ and if
$|\l|+|\m|=h(n,d)$, then
$W_\l(p)\cap \widetilde{W}_\m(p')=\emptyset$ for a general translate
$\widetilde{W}_\m(p')$ of $W_\m(p')$.
\end{cor}

\noindent
As in Bertram's paper \cite{Ber}, we cover $LQ_d$ by images of Grassmann
bundles over smaller Quot schemes.

\begin{defn}
We let $\pi_c\colon G_c\to \bP^1\times LQ_{d-c}$
denote the Grassmann bundle of $c$-dimensional
quotients of the universal bundle $\E$ on $\bP^1\times LQ_{d-c}$.
The morphism $\beta_c\colon G_c\to LQ_d$ is given as follows.
On $G_c$ there is the universal quotient bundle $\F_c$, of rank $c$.
If $i_c$ denotes the morphism $G_c\to \bP^1\times G_c$ given by
$({\rm pr}_1\circ\pi_c,{\rm id})$, then there is a natural
morphism of sheaves
$({\rm id}\times ({\rm pr}_2\circ \pi_c))^*\E\to i_{c{*}}\pi_c^*\E$
on $\bP^1\times G_c$.
We define $\E_c$ to be the kernel of this morphism composed with
$i_{c{*}}$ of the morphism to the universal quotient bundle on $G_c$:
$$\E_c=\ker(({\rm id}\times ({\rm pr}_2\circ \pi_c))^*\E\to i_{c{*}}\F_c).$$
Now $\E_c$ is a subsheaf of $\cO\otimes V$ on $\bP^1\times G_c$,
with cokernel flat over $G_c$, and hence there is an induced morphism
$\beta_c\colon G_c\to LQ_d$.
\end{defn}

\noindent
It is natural to consider degeneracy loci with respect to the bundles $\E_c$.

\begin{defn}
The degeneracy loci $\WW_{c,\l}$ and $\WW_{c,\l}(p)$
($\lambda\in \D_n$ and $p\in \bP^1$)
are the following subschemes of $G_c$:
\begin{align*}
\WW_{c,\l} &= \{\,x\in G_c\,|\,
\rk(\E_c\to \cO\otimes V/F^\perp_{n+1-\l_i})_x\lequ n+1-i-\l_i,\\
&\qquad\qquad\qquad\qquad\qquad\qquad\qquad\qquad i=1,\ldots,\ell(\l)\,\}, \\
\WW_{c,\l}(p) &= \WW_{c,\l}(p) \cap \pi_c^{-1}(\{p\}\times LQ_{d-c})
\end{align*}
\end{defn}

\subsection{Boundary structure of $LQ_d$}
\label{structure}
The boundary of $LQ_d$ is made up of points at which $\QQ$ fails to be locally
free, or equivalently, where $\E\to\cO\otimes V$ drops rank at one
or more points of $\bP^1$.
The following theorem is taken more-or-less verbatim from \cite{Ber}.

\begin{thm}
\label{structurethm}
The maps $\beta_c\colon G_c\to LQ_d$ satisfy \\
{\em (i)} Given $x\in LQ_d$, if $\QQ_x$ has rank at least $n+c$ at
$p\in \bP^1$, then $x$ lies in the image of $\beta_c$. \\
{\em (ii)} The restriction of $\beta_c$ to
$\pi_c^{-1}(\bP^1\times LM_{d-c})$ is a locally closed immersion. \\
{\em (iii)} We have
$$\beta_c^{-1}(W_\l(p))=\pi_c^{-1}(\bP^1\times W_\l(p))\cup
\WW_{c,\l}(p)$$
where on the right, $W_\l(p)$ denotes the degeneracy locus in
$LQ_{d-c}$.
\end{thm}

\begin{proof}
The argument is exactly as in \cite{Ber}.
\end{proof}

\begin{proof}[Proof of the Moving Lemma]
The assertions on $LM_d$ are clear by transversality of general
translates (cf.\ \cite[III.10.8]{hartshorne}).
So, we restrict attention to the boundary; it is enough to show that
$$\pi_c^{-1}\bigl(\bP^1\times \bigcap_{i=1}^k W_{\l^i}(p_i)\bigr)
\,\,\,{\rm and}\,\,\,
\pi_c^{-1}\bigl(\bP^1\times \bigcap_{\substack{i=1\\ i\ne j}}^k
W_{\l^i}(p_i)\bigr)\cap \WW_{c,\l^j}(p_j)\,\,\,(j=1,\ldots,k)$$
have dimension at most $h(n,d)-\sum_i|\l^i|-1$.
The first of these is taken care of by induction on the degree $d$.
For the second: there are loci
$$U_c(p)\subset \pi_c^{-1}(\{p\}\times LQ_{d-c})$$
defined by the condition $\rk(\E_c\to \cO\otimes V)=n-c$, and
evaluation maps
$$f_c^p\colon U_c(p)\to IG(n-c, 2n).$$
Given $x\in LQ_d$, if $\QQ_x$ has rank $n+c$ then $x\in \beta_c(U_c(p))$,
so it suffices to prove dimension estimates for
\begin{equation}
\label{dlocusu}
\pi_c^{-1}\bigl(\bP^1\times \bigcap_{\substack{i=1\\ i\ne j}}^k
W_{\l^i}(p_i)\bigr)\cap \WW_{c,\l^j}(p_j)\cap U_c(p_j).
\end{equation}
But now $\WW_{c,\l^j}(p_j)\cap U_c(p_j)=(f_c^p)^{-1}(\X')$ for a
Schubert variety $\X'$ in $IG(n-c,2n)$ of codimension at least
$|\l^j|-c(c+1)/2$.
Now, by transversality of a general translate, the intersection
(\ref{dlocusu}) has dimension at most
$h(n,d)-1-\sum_i|\l^i|$.
\end{proof}

\section{Intersections on $LQ_d$ and the quantum Giambelli formula}
\label{itlqd}

As outlined in \S \ref{overview}, we need three kinds of rational
equivalences on $LQ_d$.
Two of the
rational equivalences we require follow from a simple
intersection-theoretic lemma, while in the remainder of this section
we develop the third kind.
The Chow group (algebraic cycles modulo rational equivalence)
of a scheme $\X$ is denoted $A_*\X$.
We employ the following notation.

\begin{defn}
Let $p$ denote a point of $\bP^1$: \\
(i) $\ev^p\colon LM_d\to LG$ is the evaluation at $p$ morphism; \\
(ii) $\tau(p)\colon LQ_d(p) \to LQ_d$ is the projection from the
{\em relative Lagrangian Grassmannian}
$LQ_d(p) := LG_n(\QQ|_{\{p\}\times LQ_d})$, that is,
the closed subscheme
of the Grassmannian ${\rm Grass}_n$ of rank-$n$ quotients \cite{Groth}
of the indicated coherent sheaf,
defined by the isotropicity condition on the kernel of the
composite morphism from $\cO_{{\rm Grass}_n}\otimes V$ to the universal
quotient bundle of the relative Grassmannian; \\
(iii) $\ev(p)\colon LQ_d(p)\to LG$ is the evaluation morphism
on the relative Lagrangian Grassmannian; \\
(iv) $\ev_c^p\colon \pi_c^{-1}(\{p\}\times LM_{d-c})\to
IG(n-c,2n)$ is evaluation at $p$.
\end{defn}

\begin{lemma}
\label{simplepullbackprop}
Let $T$ be a projective variety which is a homogenous space
for an algebraic group $G$.
Let $\X$ be a scheme, equipped with an action of the group $G$.
Let $U$ be a $G$-invariant integral open subscheme of $\X$, and let
$f\colon U\to T$ be a $G$-equivariant morphism.
Then the map on algebraic cycles
\begin{equation}
\label{simplepullbackeq}
[V]\mapsto \bigl[\,f^{-1}(V)\overline{\phantom j}\,\bigr]
\end{equation}
{\em(}this sends the cycle of an integral closed subscheme $V$ of $T$ to the
cycle associated with the closure in $\X$ of the pre-image in $U${\em)}
respects rational equivalence, and hence induces a map on Chow groups
$A_*T\to A_*\X$.
\end{lemma}

\begin{proof}
Denote by $Z$ the closure in $\X\times T$ of the graph of $f$.
Since $T$ is projective, the map $Z\to \X$ is proper.
The action of $G$ on $U$ extends to an action on $Z$.
The map $Z\to T$ is $G$-equivariant, and hence is flat; the restriction
to $Z\ssm U$ is flat as well, of smaller relative dimension.
It follows that the map on cycles (\ref{simplepullbackeq})
is equal to the composite of flat pullback to $Z$, followed by
proper pushforward to $\X$.
Since these operations respect rational equivalence, so does the composite.
\end{proof}

Let us denote by $\wt{\X}_{\m}$, $\wt{\X}'_{\m}$ and 
$\wt{\X}''_{\m}$ the translates of $\X_{\m}$, $\X'_{\m}$ and
$\X''_{\m}$ by a general element of the algebraic group $Sp_{2n}$.

\begin{cor}
\label{neededratequivcor}
The following cycles are rationally equivalent to zero on $LQ_d$, and also
on $LQ_d(p')$, for any $\l\in \D_n$, with $\ell=\ell(\l)$ and
$r=2\lfloor(\ell+1)/2\rfloor$,
and distinct points
$p$, $p'\in \bP^1$.\\
{\em(i)} {\em(}$\ell\gequ 3${\em)}
$\bigl[\,(\ev^p)^{-1}(\X_\l)\overline{\phantom j}\,\bigr]-
\sum_{j=1}^{r-1} (-1)^{j-1}
\bigl[\,(\ev^p)^{-1}(\X_{\l_j,\l_r}\cap
\wt{\X}_{\l\smallsetminus\{\l_j,\l_r\}})\overline{\phantom j}\,\bigr].$ \\
{\em (ii)} {\em(}$\ell\gequ 3${\em)}
$\sum_{j=1}^{r-1} (-1)^{j-1}
\bigl[\,\beta_1\bigl((\ev_1^p)^{-1}(\X'_{\l_j,\l_r}\cap
\wt{\X}'_{\l\smallsetminus\{\l_j,\l_r\}})\bigr)\overline{\phantom j}\,\bigr].$ \\
{\em (iii)} {\em(}$r=\ell\gequ 4${\em)}
$\sum_{j=1}^{r-1} (-1)^{j-1}
\bigl[\,\beta_2\bigl((\ev_2^p)^{-1}(\X''_{\l_j,\l_r}\cap
\wt{\X}''_{\l\smallsetminus\{\l_j,\l_r\}})\bigr)\overline{\phantom j}\,\bigr].$
\end{cor}

\begin{proof}
We get (i) from the classical Giambelli formula by
invoking Lemma \ref{simplepullbackprop} with
$\X=LQ_d$ or $\X=LQ_d(p')$, and $U=LM_d$.
The rational equivalences of Corollary \ref{forcor} give us (ii) and (iii):
we invoke Lemma \ref{simplepullbackprop} with
$\X=\pi_c^{-1}(\{p\}\times LQ_{d-c})$
or $\X=\pi_c^{-1}(\{p\}\times LQ_{d-c})\times_{LQ_d}LQ_d(p')$, and
$U=\pi_c^{-1}(\{p\}\times LM_{d-c})$,
where $c=1$ for (ii) and $c=2$ for (iii).
\end{proof}

\smallskip
Recall that given a closed subscheme $Z$ of a scheme $\X$,
one defines $[Z]\in A_*\X$ to
be the class in the Chow group of the formal sum $\sum m_V [V]$ where
$V$ runs over integral closed subschemes of $\X$ which are irreducible
components of $Z$, and where $m_V$ is the geometric multiplicity of
$V$ in $Z$.
In the following statement, we write
$[Z]_k$ to denote the $k$-dimensional component of $[Z]$,
that is, the class of the 
same sum with $V$ restricted to the irreducible components
of $\X$ which have dimension $k$.

\begin{prop}
\label{pointmove}
{\em (a)} Suppose $\l$ and $\m$ are in $\D_n$, and let $p$, $p'$, $p''$ be
distinct points in $\bP^1$.
Assume $\ell(\l)$ equals $1$ or $2$
and $\ell(\m)\gequ 2$.
Let $k=h(n,d)-|\l|-|\m|$.
Then
\begin{align}
\bigl[W_\l(p)\cap \wt{W}_\m(p')\bigr]_k
&= \bigl[W_\l(p)\cap \wt{W}_\m(p) \bigr]_k    \ \ \   {\it in}\ A_*LQ_d,
                                                  \label{pppone}    \\
\bigl[\tau(p'')^{-1}\bigl(W_\l(p)\cap \wt{W}_\m(p')\bigr)\bigr]_k
&= \bigl[\tau(p'')^{-1}\bigl(W_\l(p)\cap \wt{W}_\m(p) \bigr)\bigr]_k
    \ \ \   {\it in}\ A_*LQ_d(p''),
                                                  \label{ppptwo}
\end{align}
where $\wt{W}_{\m}(p)$
denotes degeneracy locus with respect to a general translate of
the isotropic flag of subspaces.

\medskip
\noindent
{\em (b)} In $A_* LQ_d$, we have
\begin{align}
\bigl[W_\l(p) \cap \wt{W}_\m(p)\bigr]_{k} 
&= \bigl[\,(\ev^p)^{-1}(\X_\l\cap
\wt{\X}_\m)\overline{\phantom j}\,\bigr] \label{threecycles} \\
&\quad {}+
\bigl[\,\beta_1\bigl((\ev_1^p)^{-1}(\X'_\l\cap
\wt{\X}'_\m)\bigr)\overline{\phantom j}\,\bigr] \nonumber \\
&\quad {}+\delta_{\ell,2}
\bigl[\,\beta_2\bigl((\ev_2^p)^{-1}(\X''_\l\cap
\wt{\X}''_\m)\bigr)\overline{\phantom j}\,\bigr], \nonumber
\end{align}
(with $\delta_{\ell,2}$ the Kronecker delta),
and in $A_* LQ_d(p'')$, the cycle class
$\bigl[\tau(p'')^{-1}\bigl(W_\l(p)\cap \wt{W}_\m(p) \bigr)\bigr]_k$
is equal to the right-hand side of {\em (\ref{threecycles})}.
\end{prop}

\begin{proof}
To prove (\ref{pppone}), we consider the subscheme
$Y:=(\bP^1\times W_{\l}(p))\cap \wt{W}_{\m}$ of $\bP^1\times LQ_d$.
By generic flatness and equivariance for automorphisms of
$\bP^1$ fixing the point $p$, the morphism
$Y\to \bP^1$ is flat over $\bP^1\ssm\{p\}$.
In fact, we claim, the morphism is flat on the complement of
some closed subscheme of dimension $< k$.
Restricting our attention to the $(k+1)$-dimensional irreducible
components of $Y$, now,
the rational equivalence $\{p\}\sim\{p'\}$ on $\bP^1$
pulls back to (\ref{pppone}) on $LQ_d$.

We justify this claim and at the same time establish part (b).
First, for any $c\gequ 1$,
the intersection of two general-position Schubert varieties in
$IG(n-c,2n)$, defined by the 
rank conditions in (\ref{schubert}) corresponding
to the partitions $\l$ and $\m$, has codimension at least
$|\l|+|\m|-c(c+1)$, with equality only if $\ell(\l)\gequ c$ by
Proposition \ref{meetsmallergrass}.
It now follows, by a dimension count similar to the one in the proof
of the Moving Lemma, that the irreducible components of
$W_\l(p) \cap \wt{W}_\m(p)$ of dimension $k$
are precisely the ones indicated on the right-hand side of
(\ref{threecycles}).
A simple argument using Kontsevich's spaces of stable maps
(see, e.g., \cite{FPa})
shows that each of these components lies in the closure of
$((\bP^1\ssm\{p\})\times W_{\l}(p))\cap(\bP^1\ssm\{p\})\times_{\bP^1}
\wt{W}_{\m}$.
For instance, consider
a point $x\in \beta_1((\ev_1^p)^{-1}(\X'_\l\cap
\wt{\X}'_\m))$.
Such a point corresponds to a degree $(d-1)$ map $\psi\colon \bP^1\to LG$
with the resulting sheaf sequence modified at the point $p$.
Let $\psi(p)\in LG$ correspond to the Lagrangian subspace $\Sigma$.
The modification selects $\Sigma'\subset \Sigma$ of codimension $1$;
this corresponds to a line $E\subset LG$.
In the generic case, this line meets $\X_\l$ in a unique point $r$
and $\wt{\X}_\m$ in a unique point $s$, both distinct from $\psi(p)$.
The map from two $\bP^1$'s joined at a point to $\psi(\bP^1)\cup E$
is a point of the space of stable maps; in fact, it is some
$$\tilde x\in \ov{M}_{0,4}(LG,d)
\times_{(LG\times LG)}(\X_{\l}\times\wt{\X}_{\m}),$$
where we let marked points $1$ and $2$ land on $r$ and $s$, respectively,
and take marked points $3$ and $4$ on the other component, distinct
from $p$.
By properness and
transversality of a general translate,
$\tilde x$ is the value at the special point of
a family of such stable maps over $\Spec\C[[t]]$, which is
generically a smooth curve mapping to $LG$.
Applying the morphism which forgets the first marked point, we get a
$\C((t))$-valued point of $M_{0,3}(LG,d)=LM_d$ with limit point $x$
in the Quot scheme.

It remains only to show that the components indicated on the
right-hand side of (\ref{threecycles}) are generically reduced.
This is clear for the first-listed component, by transversality of
a general translate.
For the others, it is a linear algebra exercise to check,
from the definition using rank conditions, that the
scheme $W_{\l}(p)\cap \wt{W}_{\m}(p)$ in the neighborhood of
a point $\beta_1(x)$ for general
$x\in (\ev_1^p)^{-1}(\X'_\l\cap \wt{\X}'_\m)$ 
(resp.\ a point $\beta_2(x)$ for general
$x\in (\ev_2^p)^{-1}(\X''_\l\cap \wt{\X}''_\m)$)
is contained in the image of
the restriction of $\beta_c$ to $\pi_c^{-1}(\{p\}\times LM_{d-c})$,
where $c=1$ (resp.\ $c=2$).
Recall, by Theorem \ref{structurethm}, that the restriction of $\beta_c$
is a locally closed immersion; now, the rank conditions, near $x$ on
$\pi_c^{-1}(\{p\}\times LM_{d-c})$, are the Schubert conditions in
(\ref{schubertprime}) or (\ref{schubertprimeprime}),
and transversality of a general translate completes the argument.

This takes care of the assertions concerning cycles on $LQ_d$.
A dimension count shows that every $k$-dimensional irreducible component in
$\tau(p'')^{-1}(W_\l(p)\cap \wt{W}_\m(p))$ and in
$\tau(p'')^{-1}(W_\l(p)\cap \wt{W}_\m(p'))$ meets
the open set where $\tau(p'')$ is an isomorphism.
Now, exactly the same argument as above establishes
the equality (\ref{ppptwo}) and the case of (b) which
concerns $LQ_d(p'')$.
\end{proof}

\medskip
We now establish the rational equivalences on $LQ_d$ --- and on
$LQ_d(p'')$ --- which directly imply the quantum Giambelli formula
of Theorem \ref{lgthm}.

\begin{prop}
Fix $\l\in \D_n$ with $\ell=\ell(\l)\gequ 3$.
Set $r=2\lfloor(\ell+1)/2\rfloor$.
Let $p$, $p'$, $p''$ denote distinct points in $\bP^1$.
Then we have the following identity of cycle classes
\begin{equation}
\label{quantumgiambellieqone}
\bigl[\,(\ev^p)^{-1}(\X_\l)\overline{\phantom j}\,\bigr]=
\sum_{j=1}^{r-1} (-1)^{j-1}
\bigl[\,\bigl((\ev^p)^{-1}(\X_{\l_j,\l_r})\cap
(\ev^{p'})^{-1}(\wt{\X}_{\l\smallsetminus\{\l_j,\l_r\}})\bigr)
\overline{\phantom j}\,\bigr],
\end{equation}
both on $LQ_d$ and on $LQ_d(p'')$, where
$\wt{\X}_\m$ denotes the translate of $\X_\m$ by a generally chosen element of
the group $Sp_{2n}$.
\end{prop}

\begin{proof}
Combining parts (a) and (b) of Proposition \ref{pointmove} gives
\begin{align*}
\bigl[\,\bigl((\ev^p)^{-1} (\X_{\l_j,\l_r}) &\cap
(\ev^{p'})^{-1}(\wt{\X}_{\l\smallsetminus\{\l_j,\l_r\}})\bigr)
\overline{\phantom j}\,\bigr] \\
&= \bigl[\,(\ev^p)^{-1}(\X_{\l_j,\l_r}\cap
\wt{\X}_{\l\smallsetminus\{\l_j,\l_r\}})\overline{\phantom j}\,\bigr] \\
&\quad {}+
\bigl[\,\beta_1\bigl((\ev_1^p)^{-1}(\X'_{\l_j,\l_r}\cap
\wt{\X}'_{\l\smallsetminus\{\l_j,\l_r\}})\bigr)\overline{\phantom j}\,\bigr] \\
&\quad {}+\delta_{\ell,r}
\bigl[\,\beta_2\bigl((\ev_2^p)^{-1}(\X''_{\l_j,\l_r}\cap
\wt{\X}''_{\l\smallsetminus\{\l_j,\l_r\}})\bigr)\overline{\phantom j}\,\bigr],
\end{align*}
for each $j$, with $1\lequ j\lequ r-1$.
Now (\ref{quantumgiambellieqone}) follows by summing and applying
(i), (ii), and (in case $\ell=r$) (iii) of Corollary \ref{neededratequivcor}.
\end{proof}

\medskip

The next result is equivalent to the quantum Giambelli relation 
(\ref{lgqgiam2}).

\begin{thm}
\label{qgthm}
Suppose $\l\in \D_n$, with $\ell=\ell(\l)\gequ 3$, and
set $r=2\lfloor(\ell+1)/2\rfloor$.
Then we have the following identity in $QH^*(LG)$:
\begin{equation}
\label{quantumgiambellieqtwo}
\sigma_\l=\sum_{j=1}^{r-1} (-1)^{j-1} \sigma_{\l_j,\l_r} 
\sigma_{\l\smallsetminus\{\l_j,\l_r\}}.
\end{equation}
\end{thm}

\begin{proof}
The classical component of (\ref{quantumgiambellieqtwo}) is the
classical Pfaffian identity.
Hence, (\ref{quantumgiambellieqtwo}) is equivalent to
\begin{equation}
\label{quantumgiambellieqthree}
0=\sum_{j=1}^{r-1} (-1)^{j-1} \langle \sigma_{\l_j,\l_r},
\sigma_{\l\smallsetminus\{\l_j,\l_r\}},
\sigma_\m \rangle_d
\end{equation}
for every $d\gequ 1$ and $\m\in \D_n$ such that
$|\l|+|\m|=h(n,d)$.
But (\ref{quantumgiambellieqthree}) follows by
applying the refined cap product operation \cite[\S8.1]{F} along $\ev(p'')$ to
a general translate of $\X_\m$ in $LG$ and to both sides of
(\ref{quantumgiambellieqone}) in $LQ_d(p'')$.
On one side, we get zero, by Corollary \ref{corzero}.
On the other side, we get the required alternating sum of Gromov--Witten
numbers, by Corollary \ref{corgw}.
\end{proof}


\section{Quantum Schubert calculus}
\label{qsc}

In this section, we use Theorem \ref{lgthm} and the algebra of
$\wt{Q}$-polynomials to find explicit 
combinatorial rules that compute some of
the quantum structure constants $e_{\l\m}^{\n}(n)$ that
appear in the quantum product of two Schubert classes.

\subsection{Algebraic background}
The $\wt{Q}$-polynomials are the duals of certain
modified Hall-Littlewood polynomials. More precisely, let
$\{P_{\l}(X;t)\}$ and $\{Q_{\l}(X;t)\}$ be
the usual Hall-Littlewood polynomials, defined as in \cite[III.2]{M}. 
Let $\{Q'_{\l}(X;t)\}$ be the adjoint basis to $\{P_{\l}(X;t)\}$ 
for the standard scalar product on $\L_n[t]$; we
have $Q'_{\l}(X;t)=Q_{\l}(X/(1-t);t)$ in the notation of $\l$-rings
(see \cite{LLT}). According to \cite[Prop.\ 4.9]{PR}, we have, for 
$\l\in\E_n$, 
\[
\dis
\wt{Q}_{\l}(X)=\omega(Q'_{\l}(X;-1)),
\]
where $\omega$ 
is the duality involution of \cite[I.2]{M}.

By property (b) of \S \ref{bdefs}
there exist integers $e(\l,\m;\,\n)$ such that
\begin{equation}
\label{Qmult}
\wt{Q}_{\l}(X)\,\wt{Q}_{\mu}(X)=\sum_{\n}e(\l,\m;\,\n)\,\wt{Q}_{\n}(X).
\end{equation}
The coefficients $e(\l,\m;\,\n)$ are independent of $n$, and defined for
any partitions $\l,\m,\nu\in\E_n$. They occur in the second
author's description of Schubert calculus for the Arakelov Chow ring
of $LG$, considered as a scheme over the integers (see \cite{T2}).

There are explicit combinatorial rules (albeit, involving signs)
for generating the numbers
$e(\l,\mu;\,\n)$, which follow by specializing corresponding formulas 
for the multiplication of Hall-Littlewood polynomials 
(see \cite[Sect.\ 4]{PR}
and \cite[III.3.(3.8)]{M}). In particular one has the following 
Pieri type formula for $\l$ {\em strict} (\cite[Prop. 4.9]{PR}):
\begin{equation}
\label{Qpieri}
\wt{Q}_{\l}(X)\,\wt{Q}_k(X)=\sum_{\m} 2^{N(\l,\m)}\,\wt{Q}_{\m}(X),
\end{equation}
where the sum is over all partitions $\m\supset\l$ with
$|\m|=|\l|+k$ such that 
$\m/\l$ is a horizontal strip, and $N(\l,\m)$ is the number of 
connected components of $\m/\l$ not meeting the first column.
Recall, by property (d) of \S \ref{bdefs}, that
\[
\wt{Q}_{\l}(X)\,\wt{Q}_n(X)=\wt{Q}_{(n,\l)}(X),
\]
for any partition $\l\in\E_n$.

When $\l$, $\mu$ and $\nu$ are strict partitions, the $e(\l,\mu;\,\n)$
are classical structure constants for $LG(n,2n)$,
\[
\s_{\l}\s_{\m}=\sum_{\n\in{\D_n}} e(\l,\mu;\,\n)\,\s_{\n},
\]
and hence are nonnegative integers. For strict
$\l$, $\mu$ and $\nu$,
Stembridge \cite{St} has given a combinatorial rule for the
numbers $e(\l,\mu;\,\n)$, analogous to the Littlewood-Richardson
rule in type $A$. More precisely, we have
\begin{equation}
\label{etof}
e(\l,\m;\,\n)=2^{\ell(\l)+\ell(\m)-\ell(\n)}f(\l,\m;\,\n)
\end{equation}
where $f(\l,\m;\,\n)$ is equal to the number of
marked tableaux of weight
$\l$ on the shifted skew shape ${\mathcal S}(\n/\m)$ satisfying certain
conditions (for details, see loc.\ cit.). We also note that
the $f(\l,\m;\,\n)$ are the 
(classical) structure constants in the cup product decomposition
\[
\tau_{\l} \tau_{\m} = 
\sum_{\n\in\D_n}  f(\l,\m;\,\n)\,\tau_{\n}
\]
where $\t_{\l}$, $\t_{\m}$ and $\t_{\n}$ denote Schubert classes 
in the cohomology of the even orthogonal
Grassmannian $OG(n+1,2n+2)$ (see \cite[Sect.\ 6]{P}). 

\begin{prop}
\label{pieriprop}
For $\l$ and $\n$ in $\D_n$ and any integer $k$ with $1\lequ k\lequ n$,
we have
\[
e(\l,k;\, (n+1,\n))= 2^{\ell(\l)-\ell(\n)}e(\n,n+1-k;\, \l).
\]
\end{prop}
\begin{proof}
This follows from the combinatorial description of the structure constants
given in (\ref{Qpieri}).
\end{proof}

\subsection{Quantum multiplication}
Recall from the Introduction that for any 
$\l,\m\in\D_n$ there is a formula 
\[
\s_{\l}\cdot \s_{\m} = 
\sum  e_{\l\m}^{\n}(n)\,\s_{\n}\, q^d
\]
in $QH^*(LG(n,2n))$. In Proposition \ref{linenumbers} we showed that
the quantum structure constants in 
degree $d=1$ for $LG(n,2n)$ are $1/2$ times classical structure constants
for $LG(n+1,2n+2)$. Using Theorem \ref{lgthm},
the Pieri rule (\ref{Qpieri}) and
Proposition \ref{pieriprop}, we now get
\begin{prop}[Quantum Pieri Rule] 
\label{qupieri}
For any $\l\in\D_n$ and $k\gequ 0$ we have
\begin{equation}
\label{quantumpieri}
\s_{\l}\cdot\s_k=\sum_{\m} 2^{N(\l,\m)}\s_{\m}+\sum_{\n}
2^{N'(\nu,\l)} \s_{\nu} \, q
\end{equation}
where the first sum is classical, as in {\em (\ref{Qpieri})}, while 
the second is over all strict $\n$ contained in $\l$ with
$|\n|=|\l|+k-n-1$ such that 
$\l/\n$ is a horizontal strip, and $N'(\nu,\l)$ is one less than the 
number of connected components of $\l/\n$.
\end{prop}

For any $d,n \gequ 0$ and partition $\n$, let $(n^d,\n)$ denote the
partition 
\[
(n,n,\ldots,n,\n_1,\n_2,\ldots),
\]
where $n$ appears $d$
times before the first 
component $\n_1$ of $\n$.
Theorem \ref{lgthm} now gives

\begin{thm}
\label{qstructthm}
For any $d\gequ 0$ and strict partitions $\l,\m,\n\in\D_n$ 
with $|\n|=|\l|+|\m|-d(n+1)$,
the quantum structure constant $e_{\l\m}^{\n}(n)$
{\em(}which is also the Gromov--Witten invariant
$\langle \s_{\l}, \s_{\m}, \s_{\n'} \rangle_d${\em)} satisfies
\[
e_{\l\m}^{\n}(n)=
2^{-d}\,e(\l,\m;\,((n+1)^d,\n)).
\]
\end{thm}

\begin{cor}
\label{positivity}
For any $d\gequ 0$ and
$\l,\m,\n\in\D_{n-1}$,
the coefficient $e(\l,\m;\,(n^d,\n))$ is a nonnegative integer divisible
by $2^d$. 
\end{cor}

\noindent
{\bf Remark.} It is not true that the coefficients $e(\l,\m;\,\n)$
in (\ref{Qmult}) are all nonnegative. For example, we have
\begin{align*}
\wt{Q}_{3,2,1}(X_4)\cdot\wt{Q}_{3,2,1}(X_4) &= 8\,\wt{Q}_{4,4,4} +
 4\,\wt{Q}_{4,3,2,2,1} + 4\,\wt{Q}_{4,2,2,2,2}-
4\,\wt{Q}_{4,4,2,2}\\
 &\qquad{}+\wt{Q}_{3,3}(-4e_1e_2e_3+2e_2^3+4e_3^2) \\
 &\qquad{}+\wt{Q}_{1,1}(4e_1^2e_4^2-4e_1e_2e_3e_4+e_2^2e_3^2).
 \end{align*}

\medskip
\noin
{\bf Example.} It follows from the Remark that
\[
\s_{3,2,1}\cdot\s_{3,2,1} = q^3
\]
in $QH^*(LG(3,6))$. This in turn implies that there is a single rational cubic 
curve passing through $3$ general points on $LG(3,6)$.
(See Corollary \ref{someproducts} for a generalization.)

\subsection{Eight-fold symmetry}
\label{8fold}
For any partition $\l\in\D_n$ with $\ell(\l)=r$, let $\l^*=
(n+1-\l_r,\ldots,n+1-\l_1)$. The next result illustrates a symmetry
enjoyed by the Gromov--Witten invariants for $LG(n,2n)$.
\begin{thm}
\label{symmetrythm}
For any $d,e\gequ 0$ with $d+e=\ell(\l)$, we have
\begin{equation}
\label{symmetry}
2^{n+d} \,\langle \s_{\l}, \s_{\m}, \s_{\n} \rangle_d
     =  2^{ \ell(\mu)+\ell(\nu)+e }\, \langle \s_{\l^*}, 
\s_{\mu'}, \s_{\nu'}\rangle_e.
\end{equation}
Also, for every $d>\ell(\l)$, we have 
$\langle \s_{\l}, \s_{\m}, \s_{\n} \rangle_d=0$.
\end{thm}

\begin{proof}
We can pass to real coefficients, and introduce
\[
\wt{\sigma}_{\l}=2^{n/4-\ell(\l)/2}\,\sigma_{\l}
\ \ \ \ \ \mathrm{and}
\ \ \ \ \ 
\wt{q} = q/2.
\]
Equation (\ref{symmetry}) is equivalent to the assertion that
the coefficient of $\wt{\sigma}_{\n'}{\wt{q}}^d$ in
$\wt{\s}_{\m}\s_{\l}$ is equal to the
coefficient of $\wt{\s}_{\m}{\wt{q}}^e$ in
$\wt{\s}_{\n'}\s_{\l^*}$.
We observe (a restatement of quantum Pieri) that
with respect to the basis
$\{\wt{\s}_\l\}$ for $H^*(LQ,\R)$, the matrices representing the operators
\begin{equation}
\label{opone}
{-}\cup\sigma_i\colon H^{2{*}}(LG,\R)\to H^{2({*}+i)}(LG,\R)
\end{equation}
and
\begin{equation}
\label{optwo}
\wt{q}\text{-coefficient}({-}*\sigma_{n+1-i})\colon
H^{2({*}+i)}(LG,\R)\to H^{2{*}}(LG,\R)
\end{equation}
(where $*$ in (\ref{optwo}) denotes the product in $QH^*$)
are transposes of each other.
Using the quantum Giambelli formulas
(\ref{lgqgiam1}) and (\ref{lgqgiam2}),
we may express the operations ${\wt{q}}^d\text{-coefficient}({-}*\sigma_{\l})$
and ${\wt{q}}^e\text{-coefficient}({-}*\sigma_{\l^*})$ as expressions
in (commuting) operators (\ref{opone}) and (\ref{optwo}), and
again they are transpose to each other. This establishes (\ref{symmetry});
the assertion that
$\langle \s_{\l}, \s_{\m}, \s_{\n} \rangle_d=0$ when $d>\ell(\l)$
is also clear.
\end{proof}

\medskip
As a consequence of the previous results, 
we obtain some cases of a quantum Littlewood-Richardson rule in type
$C$, specifically, expressions for $\s_{\l}\s_{\m}$ whenever
$\ell(\m)\lequ 3$.

\begin{cor}
\label{QLRcor}
For any $\l\in\D_n$ and integers $a,b,c$ with $a>b>c$, we have
\begin{align}
\s_{\l}\cdot\s_{a,b}=
\sum_{\m}e(\l,(a,b);\,\m)\,\s_{\m}
&+ \frac{1}{2}\sum_{\n}e(\l,(a,b);\, (n+1,\n))\,\s_{\n}\, q \label{lr2} \\
&+ \sum_{\rho} f(\rho,(a,b)^*;\,\l)\,\s_{\rho}\, q^2, \nonumber \\
\s_{\l}\cdot\s_{a,b,c}=
\sum_{\m}e(\l,(a,b,c);\,\m)\,\s_{\m}
&+ \frac{1}{2}\sum_{\n}e(\l,(a,b,c);\, (n+1,\n))\,\s_{\n}\, q \label{lr3} \\
&+ \sum_{\rho}f(\rho,(a,b,c)^*;\,(n+1,\l))\,\s_{\rho}\,q^2 \nonumber \\
&+ \sum_{\eta} f(\eta,(a,b,c)^*;\,\l)\,\s_{\eta}\, q^3. \nonumber
\end{align}
\end{cor}

In particular, we see that
all the quantum structure constants on the Lagrangian Grassmannians $LG(n,2n)$
with $n\lequ 7$ are powers of $2$ (prescribed by
Theorem \ref{symmetrythm} and Proposition \ref{linenumbers}) times
classical structure constants.

Note that in (\ref{lr2}), we could have applied
Theorem \ref{symmetrythm} to rewrite the second summand on the right
as a structure constant involving $(a,b)^*$:
\begin{equation}
\label{newid1}
e(\l,(a,b);\, (n+1,\n)) = 2^{\ell(\l)-\ell(\n)} e(\n,(a,b)^*;\,(n+1,\l)).
\end{equation}
The disappearance of the factors of $2$ upon replacing the $e$'s by
$f$'s after applying (\ref{symmetry}) is a more 
general phenomenon.
For any partitions $\l$, $\m$ and $\n$, let us {\em define} the integer
$f(\l,\m;\,\n)$ using equation (\ref{etof}). 

\begin{prop}
For any $\l$, $\m\in D_n$, we have
$$\s_{\l}\s_{\m}=\sum_{d+e=\ell(\m)}\sum_{\n\in \D_n}
f(\n,\m^*;\,((n+1)^e,\l))\,\s_{\n}\,q^d$$
in $QH^*(LG)$. In other words, any degree $d$ quantum structure constant
satisfies 
\begin{equation}
\label{newid2}
e_{\l\m}^{\n}(n)=f(\n,\m^*;\,
((n+1)^e,\l))
\end{equation}
where $d+e=\ell(\m)$.
\end{prop}

\noindent
{\bf Remarks.} a) The identities (\ref{newid1}) (which involves only
classical structure constants) and (\ref{newid2}) appear
to be new. Note that (\ref{newid2})
implies that any degree $d$ quantum structure constant
$e_{\l\m}^{\n}(n)$, where $d=\ell(\mu)$ or $d=\ell(\mu)-1$, is equal to a
classical structure constant.

\medskip
\noindent
b) For any $\l\in\E_n$ define the $\wt{P}$-polynomial
$\wt{P}_{\l}(X)=2^{-\ell(\l)}\wt{Q}_{\l}(X)$; we have
\[
\wt{P}_{\l}(X) \wt{P}_{\m}(X) = 
\sum_{\n\in\E_n}  f(\l,\m;\,\n)\,\wt{P}_{\n}(X).
\]
In a companion paper
to this one \cite{KTorth}, we show that the multiplication
of $\wt{P}$-polynomials describes the quantum cohomology ring of the
spinor variety $OG(n+1,2n+2)$. At present we do not have a combinatorial
interpretation for the coefficients
$f(\l,\m;\,(n^e,\n))$ when $e\gequ 2$.

\medskip
If we fix an integer $d$ and strict partitions 
$\l$, $\m$, $\n$, then by repeated
application of 
(\ref{symmetry}) we can relate $8$ different Gromov--Witten
invariants to each other. More precisely,
if we define operators $A$, $B$ and $C$ by
\[
A\,\langle \l,\m,\n \rangle = \langle \l^*,\m',\n'\rangle, \ \ \ \ \ 
B\,\langle \l,\m,\n \rangle = \langle \l',\m^*,\n'\rangle, \ \ \ \ \ 
C\,\langle \l,\m,\n \rangle = \langle \l',\m',\n^*\rangle
\]
then $\{A,B,C\}$ generates a $(\Z/2\Z)^3$-symmetry on the table of
Gromov--Witten invariants.
This symmetry dictates the vanishing of many of these numbers:
\begin{prop}
\label{whichnonzero}
Let $\l$, $\m$, $\n\in \D_n$ and let $d$ be an integer.
The inequalities 
\begin{equation}
\label{twoineq}
0\lequ d\lequ \ell(\l)\quad{\rm and}\quad
\ell(\l)+\ell(\m)-n\lequ d\lequ \ell(\l)+\ell(\m)+\ell(\n)-n
\end{equation}
are necessary conditions for the Gromov--Witten number 
$\langle \s_\l, \s_\m, \s_\n\rangle_d$ to be nonzero.
Moreover, if the two sides of any of the inequalities in 
{\em(\ref{twoineq})} differ by $0$ or $1$, then
$\langle \s_\l, \s_\m, \s_\n\rangle_d$ is related by
the eight-fold symmetry to a classical structure constant.
\end{prop}

As an application, we have the following rule for quantum multiplication
by $\s_{\rho_n}$:

\begin{cor}
\label{someproducts}
For any $\l\in\D_n$, we have
$\s_{\l}\cdot\s_{\rho_n}=\s_{\l^{\prime{*}}}\,q^{\ell(\l)}$
in $QH^*(LG(n,2n))$.
In particular, $\s_{\rho_n}\cdot\s_{\rho_n}=q^n$.
\end{cor}


\begin{thebibliography}{C-F2}





\bibitem[AS]{AS} A. Astashkevich and V. Sadov : 
{\em Quantum cohomology of partial flag manifolds
$F_{n_1,\ldots,n_k}$}, Comm. Math. Phys. {\bf 170} (1995), no. 3, 503--528.



\bibitem[BGG]{BGG} I. N. Bernstein, I. M. Gelfand and S. I. Gelfand :
{\em Schubert cells and cohomology of the spaces $G/P$}, Russian 
Math. Surveys {\bf 28} (1973), no. 3, 1--26. 


\bibitem[Be]{Ber} A. Bertram :
{\em Quantum Schubert calculus},
Adv. Math. {\bf 128} (1997), no. 2, 289--305.


\bibitem[BCF]{BCF} A. Bertram, I. Ciocan-Fontanine and W. Fulton :
{\em Quantum multiplication of Schur polynomials},
J. Algebra {\bf 219} (1999), no. 2, 728--746. 


\bibitem[Bo]{Bo} A. Borel : {\em Sur la cohomologie des espaces
fibr\'{e}s principaux et des espaces 
homog\`{e}nes de groupes de Lie compacts}, Ann. of Math. {\bf 57} (1953),
115--207.


\bibitem[C]{C} L. Chen :
{\em Quantum cohomology of flag manifolds},
Adv. Math. {\bf 174} (2003), 1--34.


\bibitem[C-F1]{CF1}  I. Ciocan-Fontanine :
{\em The quantum cohomology ring of flag varieties}, Trans.
Amer. Math. Soc. {\bf 351} (1999), no. 7, 2695--2729.


\bibitem[C-F2]{CF2} I. Ciocan-Fontanine : 
{\em On quantum cohomology rings of partial flag
varieties}, Duke Math. J. {\bf 98} (1999), no. 3, 485--524.


\bibitem[D1]{D1} M. Demazure :
{\em Invariants sym\'{e}triques des groupes de Weyl et torsion},
Invent. Math. {\bf 21} (1973), 287--301.


\bibitem[D2]{D2} M. Demazure :
{\em D\'esingularization des vari\'et\'es de Schubert 
g\'en\'eralis\'ees}, Ann. Scuola Norm. Sup. Pisa Cl. Sci. (4)
{\bf 7} (1974), 53--88.




\bibitem[FGP]{FGP} S. Fomin, S. Gelfand and A. Postnikov :
{\em Quantum Schubert polynomials}, J. Amer. Math. Soc. {\bf 10}
(1997), 565--596.



\bibitem[F]{F} W. Fulton :
{\em Intersection Theory}, Second edition, 
Ergebnisse der Math. {\bf 2}, Springer-Verlag, Berlin, 1998.



\bibitem[FP]{FPa} W. Fulton and R. Pandharipande :
{\em Notes on stable maps and quantum cohomology},
in Algebraic Geometry (Santa Cruz, 1995), 45--96,
Proc. Sympos. Pure Math. {\bf 62}, Part 2, Amer. Math. Soc.,
Providence, 1997.



\bibitem[FPr]{FP} W. Fulton and P. Pragacz :
{\em Schubert varieties and degeneracy loci}, 
Lecture Notes in Math. {\bf 1689}, Springer-Verlag, Berlin, 1998.


\bibitem[GK]{GK} A. Givental and B. Kim :
{\em Quantum cohomology of flag 
manifolds and Toda lattices}, Comm. Math. Phys. {\bf 168} (1995), no. 3,
609--641.



\bibitem[G1]{G} A. Grothendieck :
{\em Techniques de construction et th\'eor\`emes d'existence en
g\'eom\'etrie alg\'ebrique IV: Les sch\'emas de Hilbert},
S\'eminaire Bourbaki {\bf 13} (1960/61), no. 221.

\bibitem[G2]{Groth} A. Grothendieck :
{\em Techniques de construction en g\'eom\'etrie analytique V:
Fibr\'es vectoriels, fibr\'es projectifs,
fibr\'es en drapeaux}, in
Familles d'espaces complexes et fondements de
la g\'eom\'etrie analytique,
S\'eminaire Henri Cartan {\bf 13} (1960/61), expos\'e 12.



\bibitem[H]{hartshorne} R. Hartshorne :
{\em Algebraic Geometry}, Grad. Texts in Math. {\bf 52},
Springer-Verlag, New York, 1977.


\bibitem[HB]{HB} H. Hiller and B. Boe :
{\em Pieri formula for $SO_{2n+1}/U_n$ and $Sp_n/U_n$}, Adv. in
Math. {\bf 62} (1986), 49--67.


\bibitem[KL]{KL} G. Kempf and D. Laksov :
{\em The determinantal formula of Schubert calculus},
Acta Math. {\bf 132} (1974), 153--162.



\bibitem[K1]{K1} B. Kim :
{\em Quantum cohomology of partial flag manifolds and a
residue formula for their intersection pairings}, 
Internat. Math. Res. Notices 1995, no. 1, 1--15.

\bibitem[K2]{K2} B. Kim :
{\em On equivariant quantum cohomology}, Internat. Math. Res.
Notices 1996, no. 17, 841--851.

\bibitem[K3]{K3} B. Kim :
{\em Quantum cohomology of flag manifolds $G/B$ and quantum Toda
lattices}, Annals of Math. {\bf 149} (1999), 129--148.



\bibitem[KP]{kimpandharipande} B. Kim and R. Pandharipande :
{\em The connectedness of the moduli space of maps to
homogeneous spaces},
preprint (2000), available at math.AG/0003168.

%
%

\bibitem[KM]{KM} M. Kontsevich, Y. Manin :
{\em Gromov--Witten classes, quantum
cohomology, and enumerative geometry}, Comm. Math. Phys. {\bf 164} 
(1994), no. 3, 525--562.



\bibitem[KT1]{KT} A. Kresch and H. Tamvakis :
{\em Double Schubert polynomials and degeneracy loci for the 
classical groups}, Ann. Inst. Fourier {\bf 52} (2002), no. 6, 1681--1727.

\bibitem[KT2]{KTorth} A. Kresch and H. Tamvakis :
{\em Quantum cohomology of orthogonal Grassmannians}, 
Compositio Math., to appear.




\bibitem[LLT]{LLT} A. Lascoux, B. Leclerc and J.-Y. Thibon :
{\em Fonctions de Hall-Littlewood et polyn\^{o}mes de Kostka-Foulkes
aux racines de l'unit\'{e}}, C. R. Acad. Sci. Paris {\bf 316} (1993), 1--6.


\bibitem[LT]{LT} J. Li and G. Tian :
{\em The quantum cohomology of homogeneous varieties},
J. Algebraic Geom. {\bf 6} (1997), 269--305.



\bibitem[LP]{LP1} A. Lascoux and P. Pragacz :
{\em Operator calculus for $\wt{Q}$-polynomials and 
Schubert polynomials}, Adv. Math. {\bf 140} (1998), no. 1, 1--43.



\bibitem[M]{M} I. G. Macdonald :
{\em Symmetric Functions and Hall Polynomials}, Second edition,
Clarendon Press, Oxford, 1995.


\bibitem[P]{P} P. Pragacz :
{\em Algebro-geometric applications of Schur $S$- and $Q$-polynomials}, 
S\'{e}minare d'Alg\`{e}bre Dubreil-Malliavin 1989-1990, Lecture
Notes in Math. {\bf 1478} (1991), 130--191, Springer-Verlag, Berlin, 
1991.

\bibitem[PR]{PR} P. Pragacz and J. Ratajski :
{\em Formulas for Lagrangian and orthogonal degeneracy loci;
$\wt{Q}$-polynomial approach}, Compositio Math. {\bf 107} 
(1997), no. 1, 11--87.

\bibitem[R]{R} J. Riordan :
{\em Combinatorial Identities}, John Wiley \& Sons, New York, 1968.

\bibitem[S]{S} I. Schur :
{\em \"{U}ber die Darstellung der symmetrischen und der alternierenden
Gruppe durch gebrochene lineare Substitutionen}, J. reine angew.
Math. {\bf 139} (1911), 155--250.


\bibitem[ST]{ST} B. Siebert and G. Tian :
{\em On quantum cohomology rings of Fano manifolds and a formula of 
Vafa and Intriligator}, Asian J. Math. {\bf 1} (1997), no. 4,
679--695.


\bibitem[St]{St} J. R. Stembridge :
{\em Shifted tableaux and the projective representations of symmetric
groups}, Adv. in Math. {\bf 74} (1989), 87--134.






\bibitem[T]{T2} H. Tamvakis :
{\em Arakelov theory of the Lagrangian Grassmannian}, 
J. reine angew. Math. {\bf 516} (1999), 207-223.




\bibitem[Th]{thomsen} J. Thomsen :
{\em Irreducibility of $\overline{M}_{0,n}(G/P,\beta)$},
Internat. J. Math. {\bf 9}, no. 3 (1998), 367--376.



\bibitem[V]{V} C. Vafa :
{\em Topological mirrors and quantum rings}, Essays on mirror
manifolds, 96--119, Internat. Press, Hong Kong, 1992. 

\bibitem[W]{W} E. Witten :
{\em The Verlinde algebra and the cohomology of the
Grassmannian}, Geometry, topology, \& physics, 357--422,
Conf. Proc. Lecture Notes Geom.
Topology, IV, Internat. Press, Cambridge, MA, 1995.





\end{thebibliography}
\end{document}